\documentclass[final,onefignum,onetabnum]{siamart190516}

\usepackage{caption}
\usepackage{subcaption}
\usepackage{cases}

\usepackage{lipsum}
\usepackage{amsfonts,amssymb}
\usepackage{graphicx}
\usepackage{epstopdf}
\usepackage{algorithmic}
\usepackage{verbatim,amscd}
\usepackage{latexsym, bm, float, color}
\usepackage{mathtools}
\usepackage{tikz}
\usepackage{enumitem}

\ifpdf
  \DeclareGraphicsExtensions{.eps,.pdf,.png,.jpg}
\else
  \DeclareGraphicsExtensions{.eps}
\fi

\newsiamremark{remark}{Remark}
\newsiamremark{hypothesis}{Hypothesis}
\crefname{hypothesis}{Hypothesis}{Hypotheses}
\newsiamthm{claim}{Claim}

\newsiamthm{assumption}{Assumption}
\newsiamthm{example}{Example}
\headers{On the training and generalization of DeepONets}{S. Lee and Y. Shin}

\title{On the training and generalization of deep operator networks
	\thanks{Submitted to the editors DATE.}}

\author{Sanghyun Lee\thanks{Department of Mathematics, Florida State University, Tallahassee, FL 32306
  (\email{lee@math.fsu.edu})}
\and Yeonjong Shin\thanks{Corresponding Author. Department of Mathematics, North Carolina State University, Raleigh, NC 27695-8205 USA
  (\email{yeonjong\_shin@ncsu.edu})}
}

\usepackage{amsopn}

\DeclareMathOperator*{\argmin}{argmin}

\newcommand{\mx}{m_\text{x}}
\newcommand{\my}{m_\text{y}}
\newcommand{\ONet}{O_{\text{net}}}
\newcommand{\tONet}{\widetilde{O}_{\text{net}}}
\newcommand{\Cop}{\mathcal{C}_{\text{op}}}

\ifpdf
\hypersetup{
	pdftitle={On the training and generalization of DeepONets},
	pdfauthor={S. Lee and Y. Shin}
}
\fi
\begin{document}
\maketitle
% REQUIRED
\begin{abstract}
We present a novel training method for deep operator networks (DeepONets), one of the most popular neural network models for operators.
DeepONets are constructed by two sub-networks, namely the branch and trunk networks.
Typically, the two sub-networks are trained simultaneously,
which amounts to solving a complex optimization problem in a high dimensional space.
In addition, the nonconvex and nonlinear nature makes training very challenging.
To tackle such a challenge, we propose a two-step training method that trains the trunk network first and then sequentially trains the branch network. 
The core mechanism is motivated by the divide-and-conquer paradigm and is the decomposition of the entire complex training task into two subtasks with reduced complexity. 
Therein the Gram-Schmidt orthonormalization process is introduced which significantly improves stability and generalization ability. 
On the theoretical side, we establish a generalization error estimate in terms of the number of training data, the width of DeepONets, and the number of input and output sensors. 
Numerical examples are presented to demonstrate the effectiveness of the two-step training method, including Darcy flow in heterogeneous porous media.
\end{abstract}
	
% REQUIRED
	
% REQUIRED

\section{Introduction}
Influenced by huge empirical successes of deep learning and artificial intelligence, neural network-based computational models have been explosively developed.
The effectiveness of these approaches has been documented and reported in various research papers \cite{lu2019deeponet,lin2021operator,bi2023accurate,kurth2023fourcastnet,oommen2022learning,boster2023artificial,goswami2022deep}, with references included therein.
%Among them, operator learning has gained huge attention especially due to their potential uses in partial differential equations (PDEs).
Among them, operator learning has garnered significant attention, particularly because of its potential applications in addressing forward and backward problems involving partial differential equations (PDEs).
One feature that sets the deep learning approach aside from traditional numerical methods
is the interplay between the training process
and the inference on unseen data.
While the training could be very challenging, time-consuming, and require huge computational resources,
once it is complete,
the inference can be done 
almost instantaneously.

%The Deep Operator Networks (DeepONets) \cite{lu2019deeponet} are the first neural operator models \cite{goswami2022physics}whose general structures were motivated by the universal approximation theorem \cite{chen1995universal} of neural networks for nonlinear operators.Since then, many neural network-based models for operators have been proposed mainly to improve the performance on unseen data. Some popular models include the DeepONet variants \cite{lu2022comprehensive,wang2022improved,kissas2022learning,venturi2023svd} and the Fourier Neural Operator (FNO) variants \cite{anandkumar2020neural,li2020fourier,guibas2021adaptive,benitez2023fine}. However, the accuracy of the inferred solutions substantially depends on  many factors including but not limited to network architectures, training data, hyperparameters, and training methods.
The Deep Operator Networks (DeepONets) \cite{lu2019deeponet} represent the pioneering examples of neural operator models \cite{goswami2022physics}, with their fundamental structures drawing inspiration from the universal approximation theorem \cite{chen1995universal} for nonlinear operators. Since then, various neural network-based operator models have emerged, primarily aimed at enhancing performance on previously unseen data. Notable models encompass different versions of DeepONets \cite{lu2022comprehensive,wang2022improved,kissas2022learning,liu2021multiscale,jin2022mionet,venturi2023svd} and Fourier Neural Operator (FNO) adaptations \cite{anandkumar2020neural,li2020fourier,guibas2021adaptive,benitez2023fine,tran2023factorized}. 
Nevertheless, the accuracy of the deduced solutions significantly hinges on multiple factors, including but not limited to, network architectures, training data, hyperparameters, and training methodologies.

The current study aims to comprehend the aforementioned interplay by focusing on DeepONets, 
given its status as one of the most widely used network models for operators. This work introduces two main contributions. 
Firstly,  a novel training method is developed, enhancing both training/optimization and generalization performances. The DeepONets consist of two sub-networks, namely trunk and branch networks. Standard approaches usually train these two sub-networks monolithically as a unified entity using a first-order optimization method, such as \texttt{Adam} \cite{kingma2014adam} and other variants \cite{Ruder_16_GDoverview}.
Accordingly, the training is amount to solve a complex optimization problem defined in a very high dimensional space in addition to it being nonconvex and nonlinear.
Many difficulties encountered in practice are conjunctly related 
to training, despite the inherent expressive capabilities of DeepONets \cite{deng2021convergence,lanthaler2022error,marcati2021exponential}.
% some failure cases encountered  difficulties in learning stiff operators, characterized by factors like discontinuities,
% This challenge may partially stem from the implicit regularization effect \cite{pmlr-v80-ma18a,park2019effect,smith2021on} inherent in gradient-based optimization methods, which tend to favor smoother solutions.

To improve the conventional monolithic optimization approach,
we propose a two-step training method.
The method is motivated by the divide-and-conquer paradigm. 
The core mechanism is the decomposition of the entire complex optimization problem into two subproblems with reduced complexity.
% , wherein the training of trunk and branch networks is decoupled.
As the name implies, this method involves two steps.
The first step is devoted to learning the trunk network.
The trunk network can be thought of as the basis representing the output functions, while the branch network corresponds to appropriate coefficients. 
The first step aims at finding the basis representation through the trunk network together with the corresponding coefficients without introducing the branch network.
The coefficients found in the first step shall be the desired values for the branch networks to learn. 
Consequently, in the second step, the branch network is trained to learn the values obtained from the first step.
The Gram-Schmidt orthonormalization is applied to the trunk network in the first step, which significantly improves stability and generalization ability.
This separation of the learning process
results in each step involving an optimization task with considerably reduced complexity.
% compared to the entire optimization problem.
In addition, we show that the minimized loss achieved by the two-step 
is identical to the one by 
the monolithic approach (Theorem~\ref{thm:equi-two-mono}).
We also prove that a zero training loss can be achieved under mild overparameterization of DeepONets (Theorem~\ref{thm:zero-loss-deeponet}).
Numerical examples demonstrate that the two-step training method significantly enhances the performance of DeepONets over the same dataset when compared to results obtained through monolithic training methods.

The other contribution is the generalization error estimate. The total error of DeepONets can be decomposed by approximation,
estimation, and optimization errors.
The approximation error of DeepONets is relatively well established, e.g. see \cite{deng2021convergence,lanthaler2022error,marcati2021exponential},
while the remaining two errors remain elusive.
One reason may deviate from 
that (1) every input function has to be transformed
into a finite-dimensional representation (e.g. discretization or generalized Fourier coefficients),
and (2) DeepONets learn only from a finite number of data.
Unlike function approximation tasks,
the input space for operators is, in principle, infinite-dimensional,
which requires one to appropriately extend the space-filling argument \cite{shin2020convergence}  
to the infinite-dimensional counterpart.
This is closely related to
the input domain on which DeepONets can learn and generalize.
To mathematically characterize these features,
we introduce Assumption~\ref{assmpt-inputspace}
and present multiple examples that satisfy it.
% The trunk networks can be thought of as
% the basis functions representing the output functions,
% and the branch networks are thought of as the corresponding coefficients.
We then connect the DeepONet structures to
the least squares regression \cite{cohen2013stability,cohen2019correction}
with respect to the number of random sensor points for the output functions.
Since the two-step training allows the branch network to learn 
the optimal least-squares coefficients directly, 
an existing error bound of the least-squares regression \cite{cohen2013stability,cohen2019correction}
provides the generalization error estimate for DeepONets.
% , thanks to the least-squares error bound.
Yet, there are certain conditions to be met for the trunk network,
which we formulate on Assumption~\ref{assmpt-trunk}.

The rest of the paper is organized as follows.
Upon introducing the problem set and preliminaries in Section~\ref{sec:2},
the proposed two-step training method for DeepONets
is presented in Section~\ref{sec:3}.
A generalization error analysis is given in Section~\ref{sec:4}.
Numerical examples were presented in Section~\ref{sec:5}
before conclusions.

\section{Problem Setup and Preliminaries} \label{sec:2}
Let $\Omega_x \subset \mathbb{R}^{d_x}$
and $\Omega_y \subset \mathbb{R}^{d_y}$
be computational (compact) domains of interest.
Let $(\mathcal{X},d_{\mathcal{X}})$ be a metric space of functions
defined on $\Omega_x$ and let $(\mathcal{Y},\|\cdot\|_{\mathcal{Y}})$ 
be a normed vector space of functions on $\Omega_y$.
Let 
\begin{equation*}
    \mathcal{G}:\mathcal{X} \ni f \mapsto 
    \mathcal{G}[f] \in \mathcal{Y}
\end{equation*}
be the operator of interest
to be approximated by neural networks.

%Any input function $f \in \mathcal{X}$ has to be converted to a finite-dimensional quantity in order for it to be used as input for DeepONets.
In order to utilize any input function $f \in \mathcal{X}$ with DeepONets, it must first undergo a conversion into a finite-dimensional quantity.
There are various ways of extracting information of $f$.
One naive way is to consider a fixed set of discretization points
$\{x_i\}_{i=1}^{\mx} \subset \Omega_x$
and discretize $f \in \mathcal{X}$ by $\bm{f} = (f(x_1),\dots,f(x_{\mx})) \in \mathbb{R}^{\mx}$ \cite{lu2019deeponet}.
Another popular way is to extract finitely many (generalized) Fourier coefficients of $f$ \cite{li2020fourier}.
For example, suppose $f(x) = \sum_{i=1}^{\infty} \hat{f}_i \phi_i(x)$
where $\langle \cdot, \cdot \rangle$ is an appropriate inner product, $\{\phi_i\}$ is an associated orthogonal basis,
and $\hat{f}_i = \langle f, \phi_i\rangle$.
Then, $\bm{f} = (\hat{f}_1,\dots,\hat{f}_{\mx})$.
For the rest of the paper, 
bold font denotes a representation in a finite-dimensional space,
and $\mx$ is referred to as the number of input function sensors.
%the bold font is understood as a finite-dimensional conversion.

The goal of operator learning is to construct 
an operator network $\mathcal{G}_{\text{NN}}$
that approximates the target operator $\mathcal{G}$.
An operator network is a mapping defined by
\begin{equation} 
    \mathcal{G}_{\text{NN}}:\mathbb{R}^{\mx} \ni \bm{f} \mapsto 
    \mathcal{G}_{\text{NN}}[\bm{f}](\cdot) \in \mathcal{Y}.
\end{equation}
%While there are various ways to design operator networks, the Deep Operator Network is the first neural operator model proposed in the literature.Also, many existing architectures of operator networks may be viewed as special cases of DeepONets. Therefore, in this work, we focus on DeepONets, the description of which is provided in the next subsection.
While various methodologies for designing operator networks are available, the DeepONets \cite{lu2019deeponet,lu2022comprehensive} place as the pioneering neural operator model proposed in the literature. This study centers its attention on DeepONets, and a detailed explanation thereof is given in the following subsection.

\subsection{Deep Operator Networks}
For $L \in \mathbb{N}$ and $\vec{\bm{n}}=(n_0,n_1,\dots,n_L) \in \mathbb{N}^{L+1}$,
a $L$-layer feed-forward neural network (NN) is 
the mapping $\mathbb{R}^{n_0} \ni x \mapsto 
    z^L \in \mathbb{R}^{n_L}$ where $z^L$ is defined recursively by
\begin{align*}
    % N:\mathbb{R}^{n_0} \ni x \mapsto 
    % z^L \in \mathbb{R}^{n_L},
    % \qquad
    z^\ell = W^\ell \sigma(z^{\ell-1}) + b^\ell, \quad 2 \le \ell \le L,
    \quad z^1 = W^1x + b^1.
\end{align*}
Here, 
$\sigma$ is a non-linear activation function that applies element-wise,
$W^\ell \in \mathbb{R}^{n_{\ell}\times n_{\ell-1}}$ 
and
$b^\ell \in \mathbb{R}^{n_\ell}$ are the weight matrix and the bias vector of the $\ell$-th layer, respectively.
The vector $\vec{\bm{n}}$ is referred to as the network architecture
and 
$\{W^\ell,b^\ell\}_{\ell=1}^L$ is to as the network parameters.
% We note that for a fixed architecture $\vec{\bm{n}}$, 
% the choice of network parameter determines the network.

% ({\color{red}Lee: no need to use $\theta$ here since later we use $\mu$ also})

%The framework of DeepONet \cite{lu2019deeponet} combines two different NNs to approximate the target operator, which are known as \textit{branch} and \textit{trunk} networks.Depending on how the branch network is constructed, DeepONet yieldsunstacked and stacked versions \cite{lu2019deeponet}. Since the stacked DeepONet is computationally expensive, this work focuses only on the unstacked DeepONet and refers it to DeepONet for the rest of the paper.
The DeepONet framework \cite{lu2019deeponet} integrates two distinct NNs to effectively approximate the target operator, referred to as the \textit{branch} and \textit{trunk} networks. The construction of the branch network leads to two variants of DeepONet: unstacked and stacked versions \cite{lu2019deeponet}. 
While the presented materials can be applicable to both versions, for the sake of simplicity, this paper concentrates on the unstacked variant, 
denoting it as DeepONet (with a slight abuse of notation) throughout the subsequent sections.
Throughout this work, $N$ shall be referred to as the \textit{width} of DeepONet.

The branch network
is a vector-valued $L_b$-layer NN
\begin{equation} \label{def:branch}
    \bm{c}(\cdot;\theta) = (c_0(\cdot;\theta),\dots,c_{N}(\cdot;\theta))^\top,
\end{equation}
whose architecture is $\vec{\bm{n}}_b = (\mx,n_1^{(b)},\dots,n_{L_b-1}^{(b)},{N+1})$
and $\theta$ represents the network parameters.
The trunk network
is a vector-valued $L_t$-layer NN defined on $\Omega_y \subset \mathbb{R}^{d_y}$
\begin{equation} \label{def:trunk}
    \bm{\phi}(\cdot;\mu)=(1,\bm{\phi}_0(\cdot;\mu))^\top,
\end{equation}
where $\bm{\phi}_0(\cdot;\mu)=(\phi_1(\cdot;\mu),\dots,\phi_N(\cdot;\mu))$
is a $L_t$-layer NN
whose architecture is $\vec{\bm{n}}_t = (d_y,n_1^{(t)},\dots,n_{L_t-1}^{(t)},N)$
and $\mu$ represents the network parameters.
% whose network  parameters are denoted by $\mu$.   

%{\color{red} Lee: $\mu$ a good notation? we define $\theta$ as the network parameter above, so maybe $\theta_\mu$? what do you think}

% whose architecture is $\vec{\bm{n}}_t = (d_y,\dots,N+1)$.
% On the other hand, let $\phi_j(\cdot;\mu)$, $j=1,\dots,N$ be functions defined on $\Omega_y \subset \mathbb{R}^{d_y}$ by trunk network,
% that could potentially be determined by network parameter $\mu_j$.
% Several specific choices of $\{\phi_j\}$ are discussed shortly. 
Then, a DeepONet is defined 
as the inner product of the branch and the trunk networks, i.e.,
% \begin{equation} 
%     \ONet[\cdot;\Theta]:\mathcal{X} \ni f \mapsto \bm{f} \mapsto 
%     \langle \bm{c}(\bm{f};\theta), \bm{\phi}(\cdot;\mu) \rangle 
%     % = 
%     % c_0(\bm{f};\theta) +   
%     % \sum_{j=1}^N c_j(\bm{f};\theta)\phi_j(\cdot;\mu) 
%     \in \mathcal{Y},
% \end{equation}
\begin{equation} \label{def:ONet}
    \ONet[\bm{f};\Theta](y) := \bm{\phi}^\top(y;\mu)\bm{c}(\bm{f};\theta)
    = c_0(\bm{f};\theta) + 
    \sum_{j=1}^N c_j(\bm{f};\theta)\phi_j(y;\mu),
\end{equation}
where $\Theta = \{\mu, \theta\}$ is the set of trainable DeepONet parameters.

% An illustration of these networks are provided in Figure \ref{fig:setup} with more details. 
% {\color{blue} YJ: Can you please update Fig 1?}
% \begin{figure}[!h]
% \centering
% \begin{subfigure}[b]{0.45\textwidth}
%          \centering
%          \includegraphics[width=\textwidth]{Figures/branch.jpg}
%          \caption{Branch Network}
%          \label{fig:branch}
% \end{subfigure}
% \hspace{0.2in}
% \begin{subfigure}[b]{0.25\textwidth}
%          \centering
%          \includegraphics[width=\textwidth]{Figures/trunk.jpg}
%          \caption{Trunk Network}
%          \label{fig:trunk}
% \end{subfigure}
% \caption{Setup. {\color{red} This is just a placeholder for future modification} }
% \label{fig:setup}
% \end{figure}

% \begin{remark}
%     The trunk networks can be chosen in several ways.
%     They can be neural networks, which require the network parameter $\mu$
%     to be chosen appropriately. 
%     Or they can be either orthogonal polynomials or POD discrete basis, which means 
%     that there are no parameters for trunk networks.
% \end{remark}

% In this work, we only focus on feed-forward neural networks.
\subsection{Training of Deep Operator Networks}
% Let $\mathcal{G}$ be the operator of interest.
Let $\{f_k\}_{k=1}^K$ be a set of input functions from $\mathcal{X}$
and $u_k(\cdot) = \mathcal{G}[f_k](\cdot)$ be the corresponding output functions in $\mathcal{Y}$.
Ideally, one wishes to minimize 
\begin{align*}
    \mathcal{L}_{\text{ideal}}(\Theta) = \frac{1}{K} \sum_{k=1}^K \|\ONet[\bm{f}_k;\Theta](\cdot) - u_k(\cdot)\|_{\mathcal{Y}}^p,
\end{align*}
where $\|\cdot\|_{\mathcal{Y}}$ is the norm in $\mathcal{Y}$, and 
$p$ is a positive number that may depend on $\|\cdot\|_{\mathcal{Y}}$.
For example, if $\|\cdot\|_{\mathcal{Y}}$ is the $L_2$-norm, $p=2$.
In practice, however, the norm $\|\cdot\|_{\mathcal{Y}}$ has to be discretized,
and let $\|\cdot\|_{\mathcal{Y}_{\my}}$ be a discretized norm.
% Therefore,
% For training, the output functions are also discretized in a similar manner.
We then seek to find parameters of the DeepONet that minimizes 
the training loss $\mathcal{L}$ defined by
\begin{equation} \label{def:loss}
    \mathcal{L}(\Theta) = \frac{1}{K} \sum_{k=1}^K \|\ONet[\bm{f}_k;\Theta](\cdot) - u_k(\cdot)\|_{\mathcal{Y}_{\my}}^p.
\end{equation}
%Typically, first-order optimization methods (e.g. stochastic gradient descent and its variants) are employed to solve the optimization problem.Due to the non-convexity and non-linearity of the loss, however,they often fail to minimize the loss properly,which causes unsatisfactory performance for applications. Yet, as shown in several existing works \cite{chen1995universal,lu2019deeponet,deng2021convergence,marcati2021exponential},DeepONets are capable of approximating many nonlinear operators including those that frequently arise in physical and engineering problems, especially related to PDEs.
Typically, the optimization problem is tackled using first-order optimization methods such as stochastic gradient descent and its variants \cite{Ruder_16_GDoverview}. However, the nonconvex and nonlinear nature of the loss frequently obstructs these methods from achieving satisfactory loss minimization, leading to failures in some applications. However, as evidenced by multiple extant studies \cite{chen1995universal, lu2019deeponet, deng2021convergence, marcati2021exponential}, DeepONets have the capacity to approximate a multitude of nonlinear operators, including those commonly encountered in physical and engineering challenges, particularly those involving partial differential equations (PDEs). 
Such an expressivity of DeepONets, however, becomes null if an effective training mechanism is not available.
% tailored for the DeepONets.
% , it is crucial to have a robust training method 

% Therefore, the problem of attaining DeepONets that exhibit poor performance is mainly attributed to the lack of
% robust training methods.

\section{Method} \label{sec:3}
To address the aforementioned challenge, we present a novel training method for DeepONets.
% via {\color{blue} norm-minimization.}
For ease of discussion, we confine ourselves to the case of $\mathcal{Y} = L^p_\omega(\Omega_y)$
whose norm is defined by 
\begin{align*}
    \|g\|_{\mathcal{Y}} = \left(\int_{\Omega_y} |g(y)|^p d\omega(y)\right)^{\frac{1}{p}}, \qquad \forall g \in \mathcal{Y},
\end{align*}
where $\omega$ is a probability measure satisfying $\int_{\Omega_y} d\omega(y) = 1$,
and consider the corresponding discrete norm by Monte-Carlo sampling
\begin{align*}
    \|g\|_{\mathcal{Y}_{\my}} = \left(\frac{1}{\my} \sum_{i=1}^{\my} |g(y_i)|^p\right)^{\frac{1}{p}}, \qquad \forall g \in \mathcal{Y},
\end{align*}
where $\{y_i\}_{i=1}^{\my}$ are independent and identically distributed  (i.i.d.) random samples from $\omega$.
Here $\my$ is referred to as the number of output function sensors. 

For any $g \in \mathcal{Y}$, let $\bm{g} = (g(y_1),\dots,g(y_{\my}))^\top$
be the discretization of $g$.
Note that 
the discrete norm $\|g\|_{\mathcal{Y}_{\my}}$ is the standard $\ell_p$-norm $\|\bm{g}\|_{\ell_p}$ upto a multiplicative constant.
Then the set of training data can be written as
\begin{align*}
    (\bm{f}_k, \bm{u}_k) = (f_k(x_1),\dots,f_k(x_{\mx}), u_k(y_1),\dots,u_k(y_{\my})), \quad k=1,\dots, K.
\end{align*}
It follows from \eqref{def:ONet} that 
the training loss function \eqref{def:loss} is 
\begin{equation} \label{def:loss-lp}
    \begin{split}
        \mathcal{L}(\{\mu,\theta\}) 
        % &= \frac{1}{K}\sum_{k=1}^{K} \frac{1}{\my}\sum_{i=1}^{\my} 
    % |\ONet[\bm{f}_k;\Theta](y_i) - u_k(y_i)|^p \\
    &= \frac{1}{K}\sum_{k=1}^{K} \frac{1}{\my}\sum_{i=1}^{\my} 
    |\bm{\phi}^\top(y_i;\mu)\bm{c}(\bm{f}_k;\theta)  - u_k(y_i)|^p,
    \end{split}
\end{equation}
where 
${\theta}$ and $\mu$ represent all the network parameters of branch $\bm{c}(\cdot;\theta)$ and trunk $\bm{\phi}(\cdot;\mu)$ networks, respectively.

\subsection{Loss Function via Matrix Representation}

By exploiting matrix representation,
the loss function \eqref{def:loss-lp}
can be written in a simple matrix form.
Let 
% $\bm{C}({\theta}) = [\bm{c}(\bm{f}_1;\theta),\dots,\bm{c}(\bm{f}_K;\theta)] \in \mathbb{R}^{(N+1) \times K}$,
\begin{align*}
    \bm{\Phi}({\mu}) = 
    \begin{bmatrix}
        \bm{\phi}^\top(y_1;\mu) \\ \vdots \\
        \bm{\phi}^\top(y_{\my};\mu)
    \end{bmatrix}
    % = 
    % \begin{pmatrix}
    % 1 &\phi_1(y_1;{\mu}) & \cdots & \phi_N(y_1;{\mu}) \\
    % \vdots & \vdots & \ddots & \vdots \\
    % 1 & \phi_1(y_{\my};{\mu}) & \cdots & \phi_N(y_{\my};{\mu})
    % \end{pmatrix} 
    \in \mathbb{R}^{\my \times (N+1)}, \hspace{0.1cm}
    \bm{C}({\theta}) = [\bm{c}(\bm{f}_1;\theta),\dots,\bm{c}(\bm{f}_K;\theta)] \in \mathbb{R}^{(N+1) \times K},
\end{align*}
and $\bm{U} = [\bm{u}_1,\dots,\bm{u}_K] \in \mathbb{R}^{\my \times K}$.
The matrix $\bm{\Phi}(\mu)$ may be viewed as a Vandermonde-like matrix,
while $\bm{C}(\theta)$ may be viewed as the corresponding coefficient-like matrix.
It can be checked that the loss \eqref{def:loss-lp} can be expressed as
\begin{equation*} 
    \mathcal{L}(\{\mu,\theta\}) 
    % =\frac{1}{K\my}\sum_{k=1}^K \|\bm{\Phi}({\mu}) \bm{c}_k({\theta}) - \bm{u}_k\|_{\ell_p}^p
    = \frac{1}{K\my}\|\bm{\Phi}({\mu})\bm{C}({\theta}) - \bm{U}\|_{p,p}^p,
\end{equation*}
where $\|\cdot\|_{p,p}$ is the entry-wise matrix norm.
Hence, training DeepONets means solving the following optimization problem:
%Therefore, the training of DeepONets is to solve the optimization problem
\begin{equation} \label{eqn:DeepONet-train}
    \min_{\mu, \theta} \mathcal{L}(\{\mu,\theta\}):=\|\bm{\Phi}({\mu})\bm{C}({\theta}) - \bm{U}\|_{p,p}^p.
\end{equation}

% $\Phi_\mu$ is the matrix of size $\my \times N$ whose $(i,j)$-component is $\phi_j(y_i;\mu_j)$,
% $C = [\bm{c}_1,\dots,\bm{c}_K]$, and
% $U = [\bm{u}_1,\dots,\bm{u}_K] \in \mathbb{R}^{\my \times K}$.
% Here 

\subsection{Reparameterization of DeepONets}
A key part of the proposed training method is 
a specific (re)parameterization of DeepONets.
% Let us consider a specific reparameterization of the ONet.
Let $T$ be a trainable square matrix of size $N+1$
and consider a new trunk network $\hat{\bm{\phi}}$ which has the form of
$$\hat{\bm{\phi}}(\cdot;\mu, T) = T^\top \bm{\phi}(\cdot;\mu),$$
where $\bm{\phi}(\cdot;\mu)$ is the standard trunk network defined in \eqref{def:trunk}.
The resulting DeepONet is then
\begin{align*}
    \ONet[\bm{f}](y) = \hat{\bm{\phi}}^\top(\cdot;\mu, T)\bm{c}(\bm{f};\theta)  = \bm{\phi}^\top(\cdot;\mu) T \bm{c}(\bm{f};\theta).
\end{align*}
This can be also viewed as a DeepONet with 
the same trunk network $\bm{\phi}(\cdot;\mu)$ as before
but a new branch network $T\bm{c}(\cdot;\theta)$.
% This can be also viewed as a reparameterization.
In fact, if $\bm{c}$ is a standard vector-valued neural network
constructed by the architecture $\vec{\bm{n}}_b$,
$T\bm{c}(\bm{f};\theta)$ can be viewed as a reparameterization of the last layer's weight matrix and bias vector as
\begin{align*}
    T\bm{c}(\bm{f};\theta) = TW^{L_b}\sigma(z^{L_b-1}) + Tb^{L_b}.
\end{align*}
With this new parameterization, we consider the following loss function  
\begin{equation} \label{def:loss-lp-repara}
    \min_{\mu,T,\theta} 
    \mathcal{L}(\mu,T,\theta):=\|\bm{\Phi}({\mu}) T \bm{C}(\theta) - \bm{U}\|_{p,p}^p,
\end{equation}
and the proposed training method aims at solving \eqref{def:loss-lp-repara}.
We note that $\bm{\Phi}(\mu)T$ can be viewed as a Vandermonde-like matrix constructed from the new trunk network $\hat{\bm{\phi}}$ as basis.

% and the proposed training method aims at solving \eqref{def:loss-lp-repara}.
% For the rest of the paper, the above loss function is referred to as 
% $\mathcal{L}(\mu,T,\theta)$.

\subsection{Proposed Two-Step Training Method} \label{subsec:2ST}
Assuming $\my > N$,
we propose a two-step training method for solving \eqref{def:loss-lp-repara}.

\textbf{Step 1.} The first step trains the new trunk network $\hat{\bm{\phi}}$ through
the following minimization problem:
\begin{equation} \label{train-trunk}
    \min_{{\mu}, A} \mathcal{L}(\mu, A):= \|\bm{\Phi}({\mu}) A - \bm{U} \|_{p,p}^p
    \quad \text{where} \quad A \in \mathbb{R}^{(N+1)\times K}.
\end{equation}
Let $({\mu}^*, A^*)$ be an optimal solution
and let $\bm{\Phi}({\mu}^*)$ be full rank.
We then set $T^* = (R^*)^{-1}$ as the inverse of $R^*$ obtained from 
a QR-factorization of $\bm{\Phi}({\mu}^*)$, i.e., $Q^*R^* = \bm{\Phi}({\mu}^*)$.
The trunk network is then fully determined as $\hat{\bm{\phi}}(\cdot;\mu^*, T^*)$.

\textbf{Step 2.} 
The second step trains the branch network to fit $R^*A^*$.
Specifically,
we consider the optimization problem of
\begin{equation} \label{train-branch}
    \min_{\theta}  \| \bm{C}(\theta) - R^*A^* \|_{2,2}^2.
\end{equation}
%Let $\theta^*$ be an optional solution.The fully trained branch network is then $\bm{c}(\cdot;\theta^*)$.
Assuming $\theta^*$ to be an optimal solution, the fully trained branch network is given by $\bm{c}(\cdot;\theta^*)$.

\begin{remark}
    The first step replaces the use of the branch network from \eqref{def:loss-lp} to the corresponding value matrix $A$.
    The trunk network's parameters $\mu$ and the value matrix $A$ are then trained simultaneously to minimize the loss \eqref{train-trunk}.
    Since the trunk loss function is convex with respect to $A$ (assuming $p \ge 1$), 
    the first step shall avoid any difficulties caused by nonlinearity and nonconvexity from the branch network.
    In addition, the number of parameters involved with it is $|\mu| + |A|$,
    while the one for the standard loss \eqref{def:loss-lp} is $|\mu|+|\theta|$.
    Hence, as long as $|A| = (N+1)K < |\theta|$, the first step yields an optimization problem whose dimension is 
    smaller than the one of the standard loss.
    Altogether, the first step \eqref{train-trunk} is designed to facilitate training of the trunk network.
\end{remark}

\begin{remark}
    The role of $T^*$ may be viewed as 
    applying the Gram–Schmidt orthonormalization process
    on the standard trunk network 
    with respect to the discretization points $\{y_j\}$.
    The approximation capability remains the same with or without $T^*$, 
    however,
    we found numerically that the introduction of $T^*$
    significantly improves stability and generalization ability.
    % The resulting optimal trunk network is 
    % then given by 
    % $\hat{\bm{\phi}}(\cdot;\mu^*,T^*) = (R^*)^{-\top}\bm{\phi}(\cdot;\mu^*)$.
    % For a given set of new points $\{\tilde{y}_j\}_{j=1}^m$,
    % the computation of the trunk network is then carried out 
    % by solving the linear systems of equations
    % \begin{align*}
    %     (R^*)^\top [\hat{\bm{\phi}}(\tilde{y}_1;\mu^*,T^*), \dots, \hat{\bm{\phi}}(\tilde{y}_m;\mu^*,T^* )] = 
    %     [\bm{\phi}(\tilde{y}_1;\mu^*), \dots, \bm{\phi}(\tilde{y}_m;\mu^*)],
    % \end{align*}
    % and the desired trunk network evaluations 
    % are obtained.
\end{remark}

\begin{remark}
    %The proposed method  splits the entire optimization problem into  the two smaller problems, whose (computational and/or optimization) complexities are significantly reduced. In the numerical test section, we demonstrate  the advantage of our proposed method by comparing the results of a monolithic approach  that trains both networks simultaneously.
    The proposed method splits the entire optimization problem into two smaller problems, the complexities of which (either computational and/or optimization) are significantly reduced. In the numerical test section, we demonstrate the advantage of our proposed method by comparing the results of a monolithic approach that trains both networks simultaneously.
\end{remark}

% \begin{remark}
%     One of the main advantages of the above 
%     {\color{blue}a norm-based minimization} approach is that the extension of the system due to any additional training data is simple to implement. In particular, with given larger $K$, the structure of $\bm{\Phi}(\mu)$ will be the same and only the number of columns is increased in 
%     $\bm{C}({\theta})$.
% \end{remark}

\subsection{Optimization Error Analysis}
We first show that 
the loss of \eqref{def:loss-lp-repara} at an optimal solution
from the two-step training via \eqref{train-trunk} and \eqref{train-branch}
is equal to the minimum loss of the original problem \eqref{eqn:DeepONet-train}.

\begin{theorem} \label{thm:equi-two-mono}
    Suppose that the branch network's architecture is sufficiently large enough so that for any $M \in \mathbb{R}^{(N+1)\times K}$, there exists $\tilde{\theta}$ such that $\bm{C}(\tilde{\theta}) = M$.
    Let $(\mu^*,A^*)$ and $\theta^*$ be optimal solutions of \eqref{train-trunk} and \eqref{train-branch}, respectively.
    Then, 
    % $(\mu^*,\theta^*)$ minimizes the loss $\mathcal{L}$ \eqref{eqn:DeepONet-train}, i.e.,
    \begin{align*}
        \mathcal{L}(\mu^*,T^*, \theta^*) = \min_{{\mu},{\theta}} \mathcal{L}(\{{\mu},\theta\}).
    \end{align*}
\end{theorem}
\begin{proof}
    We prove the theorem by contradiction. 
    Suppose that there exists $\{\hat{{\mu}},\hat{{\theta}}\}$ that gives 
    $\mathcal{L}(\{\hat{{\mu}},\hat{{\theta}}\}) < \mathcal{L}({\mu}^*,T^*,{\theta}^*)$.
    By letting $\hat{A}:= \bm{C}(\hat{{\theta}})$, it can be checked that 
    \begin{align*}
        \mathcal{L}({\mu}^*,T^*,{\theta}^*) &= 
         \|\bm{\Phi}({\mu}^*) T^* \bm{C}({\theta}^*) - \bm{U} \|_{p,p}^p
        = \|\bm{\Phi}({\mu}^*) A^* - \bm{U} \|_{p,p}^p \\
        &\le \|\bm{\Phi}(\hat{{\mu}})\hat{A} - \bm{U} \|_{p,p}^p
        = \|\bm{\Phi}(\hat{{\mu}})\bm{C}(\hat{{\theta}})- \bm{U} \|_{p,p}^p
        = \mathcal{L}(\{\hat{{\mu}},\hat{{\theta}}\}),
    \end{align*}
    which is a contradiction.
    This shows that $\mathcal{L}(\mu^*,T^*, \theta^*) \le \min_{{\mu},{\theta}} \mathcal{L}(\{{\mu},\theta\})$. 
    
    Let $\hat{\mu} = \mu^*$
    and $\hat{\theta}$ be the branch parameter that satisfies 
    $\bm{C}(\hat{\theta}) =A^*$.
    It then can be checked that 
    \begin{align*}
        \mathcal{L}(\{\hat{{\mu}},\hat{{\theta}}\})
        = \|\bm{\Phi}({\mu}^*) \bm{C}(\hat{\theta}) - \bm{U} \|_{p,p}^p 
        = \|\bm{\Phi}({\mu}^*) A^* - \bm{U} \|_{p,p}^p 
        = \mathcal{L}({\mu}^*,T^*,{\theta}^*),
    \end{align*}
    which shows that $\min_{{\mu},{\theta}} \mathcal{L}(\{{\mu},\theta\}) \le \mathcal{L}(\mu^*,T^*, \theta^*)$.
    % Therefore, the proof is completed.
\end{proof}

Theorem~\ref{thm:equi-two-mono} assumes that the branch network can interpolate any given data.
This is indeed possible under mild conditions (e.g. sufficiently large width)
as already shown in \cite{oymak2019towards,zou2018stochastic,du2018gradient-DNN,du2018gradient-shallow,allen2018convergence}.
% one can find ${\theta}^*$ that satisfies $\bm{C}(\theta^*) = R^*A^*$.
% Note that the learning task of \eqref{train-branch} can be thought of as a type of vector-valued function approximation.
% Thus,
In addition, if the branch network is built based on two-layer networks,
the minimization problem \eqref{train-branch}
can be effectively solved by a recently developed training method, Active Neuron Least Squares (ANLS) \cite{ainsworth2021plateau,ainsworth2022ANLS}.

Lastly, for the training of the trunk network \eqref{train-trunk},
a zero loss can also be obtained if the architecture of the trunk network is appropriately chosen.
\begin{theorem} \label{thm:trunk-interpolation}
    Suppose that $p=2$, $\bm{U}$ has rank $r$
    and the trunk network $\bm{\phi}_0$ of \eqref{def:trunk} is
    a $(2\my +1)$-layer ReLU network whose architecture is given by
    \begin{align*}
        \vec{\bm{n}}_t = (d_y, 4, 4, \overbrace{\tilde{n}, \dots, \tilde{n}}^{(2\my-2) \text{times}},N),
        \quad
        \text{where} \quad 
        \tilde{n}  =2\min\{N,r\}+4.
    \end{align*}
    Then, 
    there exist $\mu^*$ and $A^*$ satisfying
    \begin{align*}
        \mathcal{L}(\mu^*, A^*):=\|\bm{\Phi}({\mu}^*) A^* - \bm{U} \|_{2,2}^2 
        \le \min_{\text{rank}(\bm{Z}) \le \min\{N,r\}} \| \bm{Z} - \bm{U}\|_{2,2}^2.
    \end{align*}
    In particular, 
    if $N \ge r$,
    we have $\mathcal{L}(\mu^*,A^*) = 0$.
\end{theorem}
\begin{proof}
    The proof can be found in Appendix~\ref{app:them:trunk-interpolation}.
\end{proof}

By combining Theorems~\ref{thm:equi-two-mono} and~\ref{thm:trunk-interpolation},
a zero training loss for DeepONets can be achieved under overparameterization of both trunk and branch networks.
\begin{theorem} \label{thm:zero-loss-deeponet}
    Suppose that the architecture of the trunk network is set as described in 
    Theorem~\ref{thm:trunk-interpolation} with $N\ge r$.
    Suppose also that 
    % the architecture of the branch network is sufficiently large enough to achieve a zero loss from \eqref{train-branch}.
    the branch network's architecture is sufficiently large enough so that for any $M \in \mathbb{R}^{(N+1)\times K}$, there exists $\tilde{\theta}$ such that $\bm{C}(\tilde{\theta}) = M$.
    Then,
    \begin{equation*}
        0 = \min_{\mu, \theta} \mathcal{L}(\{\mu,\theta\}):=\|\bm{\Phi}({\mu})\bm{C}({\theta}) - \bm{U}\|_{p,p}^p.
    \end{equation*}
\end{theorem}
\begin{proof}
    It follows from Theorem~\ref{thm:trunk-interpolation} that 
    there exist $\mu^*, A^*$ satisfying $\mathcal{L}(\mu^*,A^*) = 0$ as $N \ge r$.
    Since the branch network can achieve a zero loss from \eqref{train-branch},
    there exists $\theta^*$ such that $\bm{C}(\theta^*) = R^*A^*$.
    It then can be checked that 
    \begin{align*}
        0 &= \mathcal{L}(\mu^*,A^*) =\|\bm{\Phi}({\mu}^*) A^* - \bm{U} \|_{p,p}^p 
        =\|\bm{\Phi}({\mu}^*)T^*\bm{C}(\theta^*) - \bm{U} \|_{p,p}^p = \mathcal{L}({\mu}^*,T^*,{\theta}^*).
    \end{align*}
    The proof is then completed by Theorem~\ref{thm:equi-two-mono}.
\end{proof}

\section{Generalization Error Analysis} \label{sec:4}
% The mathematical tool we used is the spectral
% {\color{blue} YJ: Explain what we mean by the generalization error.}
The generalization error refers to 
a quantity that measures how well the \textit{learned} DeepONet performs on \textit{unseen} functions (data).
To be more precisely,
let $\mathcal{X}_K = \{f_1,\dots,f_K\} \subset \mathcal{X}$ 
be a set of functions
and $\{u_j:=\mathcal{G}[f_j]: j=1,\dots,K\}$
be the corresponding output functions of the operator $\mathcal{G}$ of interest,
which are all used for the training of DeepONets.
Let $\ONet$ be a fully trained DeepONet.
% by \eqref{eqn:DeepONet-train}.
For $f \in \mathcal{X} \backslash \mathcal{X}_K$,
the generalization error of the DeepONet $\ONet$
at $f$
is defined to be 
\begin{equation} \label{eqn:def:gen-err}
    \mathcal{E}_{\text{gen}}(\ONet[\bm{f}]) := 
    \| \mathcal{G}[f] - \ONet[\bm{f}]\|_{L^2_\omega(\Omega_y)}.
\end{equation}
The end goal of operator learning is to construct 
a neural operator $\ONet$ from
finitely many data that yield a small generalization error uniformly over $\mathcal{X}$.

In what follows, we present a generalization error analysis for DeepONets in terms of 
the number $K$ of training data,
the number $\mx$ of input function sensors,
the number $\my$ of output function sensors,
and the width $N$ of DeepONets.
The presented analysis is motivated by an approximation error analysis of DeepONets for the coefficient-to-solution map of elliptic second-order PDEs \cite{marcati2021exponential}
and we combine it with 
the analysis \cite{cohen2013stability,cohen2019correction} of the least squares approximations 
to incorporate the output training data.

Let us consider a class $\Cop$ of operators from $\mathcal{X}$
to $\mathcal{Y} = L^2_{\omega}(\Omega_y)$
% \begin{align*}
%     \mathcal{G}:\mathcal{X} \ni f \mapsto u \in \mathcal{Y} = L^2_{\omega}(\Omega),
% \end{align*}
which has a spectral form of 
\begin{align*}
    \mathcal{G}[f](y) = \sum_{j=0}^\infty \mathfrak{c}_j(f)\psi_j(y),
\end{align*}
where $\{\mathfrak{c}_j\}$'s are $L_j$-Lipschitz functionals in $\mathcal{X}'$
such that $\sum_{j=0}^\infty L_j^2 < \infty$,
and $\{\psi_j(\cdot)\}_j$ is an orthogonal basis for $L^2_{\omega}(\Omega_y)$ 
satisfying 
\begin{equation*}
    \langle \psi_i, \psi_j \rangle_{L^2_{\omega}} := \int_{\Omega_y} \psi_i(y)\psi_j(y) d\omega(y) = \delta_{ij}.
\end{equation*}
Here $\delta_{ij}$ is a Kronecker delta function, 
and $\omega$ is a probability measure on $\Omega_y$.
It then can be checked that every operator in $\Cop$ is Lipschitz.

\begin{proposition}
    Any operator $\mathcal{G} \in \Cop$ is Lipschitz continuous.
    We denote $L_{\mathcal{G}}$ as the Lipschitz constant of $\mathcal{G}$.
\end{proposition}
\begin{proof}
    For any $f, f' \in \mathcal{X}$, observe that 
    $\|\mathcal{G}[f] - \mathcal{G}[f']\|_{L^2_\omega}^2
        = \sum_{j=0}^\infty |\mathfrak{c}_j(f) - \mathfrak{c}_j(f')|^2
        \le \sum_{j=0}^\infty L_j^2 d_{\mathcal{X}}^2(f,f')$,
    which gives $\|\mathcal{G}[f] - \mathcal{G}[f']\|_{L^2_\omega} \le L_{\mathcal{G}}d_{\mathcal{X}}(f,f')$
    with $L_{\mathcal{G}} \le \sqrt{\sum_{j=0}^\infty L_j^2}$.
\end{proof}

Let $\mathcal{G}_N[f](y) = \sum_{j=0}^N \mathfrak{c}_j(f)\psi_j(y)$
be the best $N$-term approximation of $\mathcal{G}[f]$,
and $\mathcal{E}_N(\mathcal{G}[f]) := \|\mathcal{G}[f] -\mathcal{G}_N[f]\|_{L^2_{\omega}}$ be the corresponding best $N$-term approximation error.
In what follows, we make a couple of assumptions on $\mathcal{G}$ that guarantee the \textit{uniform boundedness} and the \textit{uniform decay rate} of the
best $N$-term approximation error,
inspired by \cite{marcati2021exponential}.

% For some technical details, the followings are assumed:
% \begin{assumption}
%     Let $\mathcal{G} \in C$. 
%     There is a constant $M > 0$ such that 
%     $|\mathcal{G}[f](y)| \le M$ for any $f \in \mathcal{X}$, and for almost every $y$ with respect to $\omega$.
% \end{assumption}
% Lastly, we introduce the truncation operator defined by $T_M(z) = \text{sign}(z)\max\{M,|z|\}$, and denote 
% $\tONet[f](y) := T_M(\ONet[f](y))$.

\begin{assumption}[Operators] \label{assmpt-operators}
    For any $\mathcal{G} \in \Cop$, the followings are assumed.
    \begin{enumerate}
        \item There is a constant $M > 0$ such that 
    $|\mathcal{G}[f](y)| \le M$ for any $f \in \mathcal{X}$, and for almost every $y$ with respect to $\omega$.
        \item 
        Let
        $\mathcal{E}_N(\mathcal{X}):= \sup_{f \in \mathcal{X}} \mathcal{E}_N(\mathcal{G}[f])$
        be the supremum of the best $N$-term approximation errors over $\mathcal{X}$.
        Assume that 
        $\mathcal{E}_N(\mathcal{X}) \le N^{-r_{\mathcal{G,X}}}$
        for some $r_{\mathcal{G,X}} > 0$ that depends on
        $\mathcal{X}$, $\mathcal{G}$ and the choice of basis $\{\psi_j\}$.
    \end{enumerate}
\end{assumption}

% {\color{blue} YJ: justify Assumption 4.2 by using \cite{marcati2021exponential},
% which is 
% Schwab's exponential convergence paper published in SINUM.}
We note that there are many operators of interest satisfying Assumption~\ref{assmpt-operators}.
For example, \cite{marcati2021exponential} considered the elliptic boundary value problem
\begin{align*}
    - \nabla \cdot (a \nabla u^a) = f,
\end{align*}
for some fixed source term $f$,
and studied approximation rates of 
the data-to-solution operator $\mathcal{G}:a \mapsto u^a$.
It was shown that the operator $\mathcal{G}$ 
satisfies Assumption~\ref{assmpt-operators}.

Since DeepONet requires one to extract finite-dimensional information from an infinite-dimensional class $\mathcal{X}$
for the input, 
in order to quantify how many input functions are needed to fill up 
the target domain,
we make the following assumptions.
% on the input function space $\mathcal{X}$.

\begin{assumption}[Input functions and Sensors] \label{assmpt-inputspace}
    The symbol $\lesssim$ is used to suppress constants that depend only on $(\mathcal{X},d_\mathcal{X})$.
    \begin{enumerate}
        \item 
        % Let $\mathcal{X}$ be a collection of functions of interest
        % such that for any $\epsilon > 0$, 
        For all but finitely many $\mx \in \mathbb{N}$, 
        there exist
        $\mx$ discretization points 
        $\{x_j\}_{j=1}^{\mx}$ in $\Omega_x$ satisfying 
        \begin{align*}
            \| \bm{f} - \bm{g} \|_{w,2} \lesssim  d_{\mathcal{X}}(f,g) + \mx^{-\alpha},
            \quad \forall f, g \in \mathcal{X},
        \end{align*}
        where $\|\cdot\|_{w,2}$ is a weighted Euclidean norm,
        and $\alpha >0$ is a constant that depends only on $\mathcal{X}$.
        \item 
        For any $K \in \mathbb{N}$,
        there exist $K$ input functions $\mathcal{X}_K:= \{f_1,\dots,f_K\}$ in $\mathcal{X}$ 
        satisfying
        for any $f \in \mathcal{X}$,
        \begin{equation}
            \min_{1\le k \le K} d_{\mathcal{X}}(f, f_k) \lesssim  K^{-s} + \mx^{-\alpha},
        \end{equation}
        for some $s > 0$ which depends only on $\mathcal{X}$.
    \end{enumerate}
\end{assumption}

Many input function spaces $\mathcal{X}$ used in the literature (e.g. \cite{lu2019deeponet,lu2022comprehensive,anandkumar2020neural})
satisfy Assumption~\ref{assmpt-inputspace}.
For the sake of clarity, we present three detailed examples.

\begin{example}
    Let 
    \begin{align*}
        \mathcal{X} = \left\{ f:[-1,1]^2 \to \mathbb{R} \mid \kappa \in [0, 1], f(x) = \kappa \text{ if } \|x\|_2 \le 1, \text{ and } 1 \text{ otherwise} \right\},
    \end{align*}
    and define a map $d_{\mathcal{X}}$ over $\mathcal{X}\times \mathcal{X}$ such that 
    $d_{\mathcal{X}}(f,g) := |f(0) - g(0)|$ for any $f, g \in \mathcal{X}$.
    It can be checked that $d_{\mathcal{X}}$ is a metric on $\mathcal{X}$.
    For a given set of points $\{x_j\}_{j=1}^{\mx}$ where there exists a point whose $\|\cdot\|_2$-norm is less than or equal to 1,
    let $w_j = 0$ if $\|x_j\| > 1$ and $w_j = \frac{1}{|\{x_i | \|x_i\| \le 1\}|}$ if $\|x_j\| \le 1$.
    For $f \in \mathcal{X}$, let $\bm{f} = (f(x_1),\dots,f(x_{\mx}))$
    and $\|\bm{f}\|_{w,2}^2 := \sum_{j=1}^{\mx} w_j  f(x_j)^2$.
    It then can be checked that 
    $\|\bm{f} - \bm{g}\|_{w,2} = d_{\mathcal{X}}(f,g)$ for all $f,g \in \mathcal{X}$,
    which satisfies Assumption~\ref{assmpt-inputspace}.1.
    
    For any $K \in \mathbb{N}$, let 
    $\mathcal{X}_K = \{f_1,\dots,f_K\}$ where $f_i \in \mathcal{X}$ and
    $f_i(0) = \frac{i}{K}$.
    It then can be checked that for any $f \in \mathcal{X}$, 
    there exists $g \in \mathcal{X}_K$
    such that $d_{\mathcal{X}}(f,g) \le K^{-1}$, which shows that Assumption~\ref{assmpt-inputspace}.2 holds.
\end{example}

\begin{example}
    Let 
    \begin{align*}
        \mathcal{X} = \{ f \in C^1([-1,1]) \mid  \|f\|_{C^1} \le 1 \},
    \end{align*}
    with $d_{\mathcal{X}}(f,g) = \|f-g\|_{L^2}$ for any $f,g \in \mathcal{X}$.
    Let $\{(x_j,w_j)\}_{j=1}^{\mx}$ be the Gauss-Legendre quadrature points and weights.
    Observe that for any $f \in \mathcal{X}$, we have 
    $\|\bm{f}\|_{w,2} = \|\Pi_{\mx} f \|_{L^2}$,
    where $\Pi_{\mx} f$ is the Lagrange interpolation of $f$.
    Therefore, for any $f, g \in \mathcal{X}$,
    \begin{align*}
        \|\bm{f}-\bm{g}\|_{w,2} &= \|\Pi_{\mx}(f-g)\|_{L^2}
        \le \|f-g\|_{L^2} + \|(f-g) - \Pi_{\mx}(f-g)\|_{L^2} \\
        &\lesssim d_{\mathcal{X}}(f,g) + \mx^{-1},
    \end{align*}
    which shows that Assumption~\ref{assmpt-inputspace}.1 holds.  
    
    For any $n \in \mathbb{N}$,
    let $K = (n+1)^{\mx}$,
    and consider 
    $\mathcal{X}_K = \{ f \in \mathcal{X} : \bm{f} \in \{-1 + \frac{i}{n}(2) : i=0,\dots,n\}^{\mx} \}$.
    It then can be checked that for any $f \in \mathcal{X}$,
    there exists $g \in \mathcal{X}_K$ such that
    \begin{align*}
        d_{\mathcal{X}}(f,g) &\le \|\Pi_{\mx}(f-g)\|_{L^2} + C\mx^{-1} = \|\bm{f}-\bm{g}\|_{w,2} + C\mx^{-1} \\
        &\le \|\bm{f}-\bm{g}\|_{\infty} + C\mx^{-1}
        \lesssim K^{-\frac{1}{\mx}} + \mx^{-1},
    \end{align*}
     which shows that Assumption~\ref{assmpt-inputspace}.2 holds. 
\end{example}

\begin{example}
    % Spectral convergence (generalized Fourier coefficients)
    Let $\mathcal{X} = \{f \in H^p_{\omega}([-1,1]): \|f\|_{H^p_\omega} \le 1 \}$
    where $H^p_{\omega}([-1,1])$
    is a weighted Sobolev space.
    For any $f \in H^p_{\omega}$,
    let $\bm{f} = (\hat{f}_0,\dots,\hat{f}_{\mx})$
    with $\hat{f}_k = \langle f, p_k\rangle_{L^2_\omega([-1,1])}$
    where $p_k$ is the orthonormal polynomial of degree $k$ with respect to  $\omega$.
    It then follows from the well-known spectral convergence \cite{hesthaven2007spectral} that
    for any $f, g \in \mathcal{X}$, 
    \begin{align*}
        \|\bm{f} - \bm{g}\|_2
        \lesssim \|f - g\|_{L^2_\omega}
        + \mx^{-p},
    \end{align*}
    implying Assumption~\ref{assmpt-inputspace}.1.

    For any $n \in \mathbb{N}$,
    let $K = (n+1)^{\mx}$,
    and consider 
    $\mathcal{X}_K = \{ f \in \mathcal{X} : \hat{f} \in 
    \{\frac{i}{n} : i=0,\dots,n \}^{\mx}$.
    Let $P_{\mx}f = \sum_{j=0}^{\mx} \hat{f}_j p_k$.
    It then can be checked that for any $f \in \mathcal{X}$,
    there exists $g \in \mathcal{X}_K$ such that
    \begin{align*}
        d_{\mathcal{X}}(f,g) &\le  
        \|P_{\mx}f - g\|_{L^2_\omega}
        + \|f - P_{\mx}f\|_{L^2_\omega} 
        \lesssim 
        \|\bm{f}-\bm{g}\|_{2} 
        + \mx^{-p}
        \\
        &\le \|\bm{f}-\bm{g}\|_{\infty} + \mx^{-p}
        \lesssim K^{-\frac{1}{\mx}} + \mx^{-p},
    \end{align*}
    implying Assumption~\ref{assmpt-inputspace}.2.
\end{example}

The following assumption is the one that draws a connection between
the proposed two-step training method and 
the generalization analysis. 
Roughly speaking, we generalize the assumption on the number of output sensors 
introduced in \cite{cohen2013stability,cohen2019correction}
for a class of trunk neural networks on which 
orthonormal basis can be formed with respect to a given measure defined on $\Omega_y$.
This assumption is crucial as it allows one to utilize a classical least square analysis
in the context of DeepONets.

\begin{assumption}[Trunk Networks and Sensors] \label{assmpt-trunk}
    Let $\emph{\text{F}}$ be a feasible set of trunk network parameters
    defined by
    % $\text{M} = \cup_{s = 1}^{\infty} \text{M}_{s}$
    % where 
    \begin{equation} \label{def:feasibleSet_trunk}
        \begin{split}
        \emph{\text{F}}
        &= \left\{ \mu \in \emph{\text{F}}_{t} :
        \sup_{y \in \Omega_y} \|\hat{\bm{\phi}}(\cdot;\mu,T_\mu)\|_{2}^2 < \infty \right\}
        \end{split}
    \end{equation}
    where
    $\emph{\text{F}}_{t}= \{ \mu :
        \exists~T_{\mu} \text{ such that }
        \hat{\bm{\phi}}(\cdot;\mu,T_\mu) \text{ forms orthonormal basis in }
        L^2_\omega(\Omega_y)
        \}$.
    % on which trunk network's parameters
    % are sought.
    Let 
    $\{y_1,\dots,y_{\my}\}$ be a set of discretization points randomly independently drawn from the probability measure $\omega$.
    For $r_{t} > 0$,
    suppose $\my$ is sufficiently large enough to satisfy
    \begin{equation} \label{assumption-eqn-trunk-sampling}
        \sup_{\mu \in \emph{\text{F}}} 
        \left( \sup_{y \in \Omega_y} \|\hat{\bm{\phi}}(\cdot;\mu,T_\mu)\|_{2}^2
        \right) 
        \le \kappa \frac{\my}{\log \my},
        \qquad \kappa := \frac{3\log (3/2) - 1}{2+2r_{t}}.
    \end{equation}
\end{assumption}

Lastly, we introduce assumptions for branch networks.
For simplicity, we confine ourselves to a two-layer neural network
of sufficiently large width  
so that one achieves a zero loss on \eqref{train-branch}.
\begin{assumption}[Branch Networks] \label{assmpt-branch}
    The following are assumed for branch networks.
    \begin{enumerate}
        \item 
        The branch network
        is a two-layer neural network whose activation function $\sigma$ is Lipschitz continuous with the Lipschitz constant $L_{\sigma}$.
        \item 
        For each $K$, there exists 
        a two-layer branch network of width $n_{K}$ 
        that achieves
        a zero loss \eqref{train-branch}.
        That is, 
        there exists $\theta^*$
        such that 
        $\bm{C}(\theta^*) = R^*A^*$.
        Specifically, let $\theta^* = \{\gamma_\ell, \beta_\ell,  w_\ell\}_{\ell=1}^{n_K}$
        where 
        $\gamma_\ell \in \mathbb{R}^{N+1}$,
        $w_\ell \in \mathbb{R}^{\mx}$,
        $\beta_\ell \in \mathbb{R}$.
        Then, for $k=1,\dots,K$,
        \begin{align*}
            \bm{c}(\bm{f}_k;\theta^*) := 
            \sum_{\ell=1}^{n_K} \gamma_\ell \sigma\left(
            \langle w_\ell, \bm{f}_k \rangle + \beta_\ell\right) = (R^*A^*)_k \in \mathbb{R}^{N+1},
        \end{align*}    
        where 
        $(R^*A^*)_k$
        is the $k$-th column of $R^*A^*$
        from \eqref{train-branch}.
        \item 
        % {\color{red} This may be proved given some boundedness on $\theta^*$.} 
        Let $L_{c}(K,N,\mx):=\sum_{\ell=1}^{n_K} \|\gamma_\ell\|_2\|w_\ell\|_2$.
        Suppose 
        $L_c(K,N,\mx)$ is uniformly bounded independent of $K$, $N$, and $\mx$,
        and denote its upper bound by 
        $\overline{L}_c$.
    \end{enumerate}
\end{assumption}

\begin{remark}
    The assumptions of Assumption~\ref{assmpt-branch}
    are mild and easily satisfied in many practical setups.
    The last assumption corresponds to the uniform boundedness of the Lipschitz constant for the branch networks,
    which is often used in the literature (e.g. \cite{shin2020convergence}) 
    to establish a convergence.
\end{remark}

% Let $f_1,\dots,f_K$ be input functions in $\mathcal{X}$ that are used for training,
% and let $\mathcal{X}_{K,\mx} = \{\bm{f}_1,\dots,\bm{f}_K\}$ be the set of the corresponding discretization vectors.
% For any $f \in \mathcal{X}$, let $\text{dist}(f,\mathcal{X}_{K,\mx}):= \min_{1\le i \le K} \|\bm{f} - \bm{f}_i\|_2$
% be the distance between the discretization of $f$ and the set $\mathcal{X}_{K,\mx}$ in the standard Euclidean norm.
We are now in a position to present the main theorem that characterizes 
the generalization error of the fully trained DeepONets
in terms of the number $K$ of training data,
the number $\mx$ of input function domain sensors,
the number $\my$ of output function domain sensors,
and the width $N$ of DeepONets.

\begin{theorem} \label{thm:main}
    Suppose Assumptions~\ref{assmpt-operators}, \ref{assmpt-inputspace}, \ref{assmpt-trunk}, and \ref{assmpt-branch} hold.
    Let $\ONet$ be the fully trained DeepONet,
    that is,
    the trunk networks are obtained from \eqref{train-trunk}
    with $p=2$ and $\mu^* \in \emph{\text{F}}$ defined as in Assumption~\ref{assmpt-trunk},
    and the branch network solves \eqref{train-branch}.
    Given a truncation operator $\mathfrak{T}_M(z) = \text{sign}(z)\max\{M,|z|\}$, let 
    $\tONet[\bm{f}](y) := \mathfrak{T}_M(\ONet[\bm{f}](y))$.
    Then, for any $f \in \mathcal{X}$,
    \begin{equation}
\mathbb{E}\left[\mathcal{E}_{\text{gen}}^2(\tONet[\bm{f}])\right] \lesssim 
        C(\my,r_t)N^{-r_{\mathcal{G,X},\mu^*}}
        + K^{-\alpha}
        + \mx^{-s} + \my^{-r_{t}},
        \label{eqn:thm}
    \end{equation}
    where $C(\my,r_t) =1+\frac{6\log (3/2) - 2}{(1+r_{t})\log \my}$
    and the expectation is taken over all random output function sensors $\{y_i\}_{i=1}^{\my}$.
    All the hidden constants are independent of $K$, $\mx, \my, N$ but may only depend on
    $M$, $L_{\mathcal{G}}$, $L_{\sigma}$, $\overline{L}_c$ and $(\mathcal{X},d_{\mathcal{X}})$.
    % Here, $L_{\sigma}$ and $L_{\mathcal{G}}$ are the Lipschitz constants of the activation function $\sigma$
    % and the operator $\mathcal{G}$, respectively,
    % and $L_c$ is defined in Assumption~\ref{assmpt-branch}.
\end{theorem}
\begin{proof}
    The proof can be found in Appendix~\ref{app:thm:main}.
\end{proof}

\begin{remark}
    The rate $r_{\mathcal{G,X},\mu^*}$ of convergence 
    with respect to the width $N$ of DeepONet
    is affected by the trunk network.
\end{remark}

% \begin{remark}
%     In practice, it has been widely observed and recognized that 
%     the architecture of branch networks plays a crucial role 
%     in determining the generalizability of operator networks \cite{lu2022comprehensive,raonic2023convolutional}.
%     The conditions we require for the branch networks are given in Assumption~\ref{assmpt-branch}
%     and the finiteness of the quantity $\overline{L}_c$ plays a key role.
%     % As a matter of fact, 
%     % Theorem~\ref{thm:main} indicates this practical observation 
%     % in a theoretical manner
%     % through the quantity $\overline{L}_c$
%     % for the case of two-layer neural networks.
% \end{remark}

% For the purpose of analysis, we focus on the two cases of $p=1$ and $p=2$.
% Let $p=2$ and $\bm{\mu}$ be fixed.
% % The optimal solution $A^*$ can be explicitly written as 
% % a function of $\bm{\mu}$
% Consider 
% \begin{equation}
%     \hat{A}[\bm{\mu}] := \Phi^\dagger[\bm{\mu}] U,
% \end{equation}
% where $\dagger$ represents the Moore–Penrose pseudo inverse.
% By substituting it into \eqref{train-trunk}, 
% one can consider 
% the task of finding the optimal parameters for trunk networks that solve 
% \begin{equation}
%     \min_{\bm{\mu}} \|(I - \Phi[\bm{\mu}]\Phi^\dagger[\bm{\mu}])U\|_{2,2}^2.
% \end{equation}
% It then can be shown that 
% \begin{proposition}
%     Let $p=2$.
%     Suppose that 
% \end{proposition}

% \clearpage
% \newpage

\section{Numerical Examples} \label{sec:5}
In this section, we present several numerical experiments to demonstrate the performance of the proposed two-step training method.
Throughout, the two-step training method as described in Section~\ref{subsec:2ST} is referred to as \texttt{2ST}
and the vanilla monolithic training method is referred to as \texttt{VAN}.
To illustrate the importance of the Gram-Schmidt orthonormalization (implemented by QR-factorization) in the first step, 
the one without it shall be referred to as \texttt{2STw/oQR}.
In all the numerical tests, 
we employ the standard unstacked DeepONet structure proposed in \cite{lu2019deeponet}.
% \shlee{(Lee:should we mention unstacked version)}

\textbf{Darcy's flow equation.}
In the following numerical examples, we consider Darcy's flow equation in a  bounded domain $\Omega = (-1,1)^2$ with Lipschitz boundary $\partial \Omega$:
\begin{equation} \label{eqn:Darcy}
    \begin{split}
         -\nabla \cdot (\alpha(p)\nabla p) &= f \ \text{ in } \Omega, \\
    p & = g \ \text{ on } \Gamma_D, \\
    -\alpha(p)\nabla p \cdot {\bf{n}} &= h  \ \text{ on } \Gamma_N
    \end{split}
\end{equation}
where 
$p: \Omega \rightarrow \mathbb{R}$ is the scalar pressure and 
$f$ is the body force.
A Dirichlet boundary condition 
% $g:= g(x,y)$ 
is imposed on $\Gamma_D$ 
% {\color{blue} where $|\Gamma_D| >0$,} 
and a Neumann boundary condition 
% $h := h(x,y)$ 
is given on $\Gamma_N = \partial \Omega \backslash \Gamma_D$ with ${\bf{n}}$ as the unit outward normal vector on $\Gamma_N$. 
% Here, $\Gamma = \Gamma_N \cup \Gamma_D$
We note that the conductivity $\alpha:=\alpha(p)$ could yield the equation to be nonlinear. 
In what follows, we will consider three different operators that arise from 
\eqref{eqn:Darcy}.
%Here, the Darcy's velocity vector $\bm{u}(x,y) : \Omega \rightarrow \mathbb{R}^d$ is then defined by $\bm{u}:=-\alpha(p) \nabla p$.

\textbf{Data generation.}
We employ the classical Lagrange continuous Galerkin linear finite element method (FEM) to generate the data. Both finite element libraries, deal.II \cite{deal} and FEniCS \cite{fenics}, were utilized.
Once the data is generated,
we split it into the training data
and the test data.
The training data is used for training of DeepONet and the test data is used to evaluate the performance of the trained DeepONet.

\textbf{Inference on unseen data.}
For a test input function $\bm{f}_{\text{test}}$,
% which does not belong to the training data set,
% {({\color{red} shall we define test, such as test = entire case $\backslash$ train})}
the DeepONet produces
an approximation to the corresponding output
function $u_{\text{test}}:=\mathcal{G}[f_{\text{test}}]$.
Let $\{y_i^{(\text{test})}\}_{i=1}^{M_\text{test}}$
be a set of points from $\Omega_y$
to be used for evaluating the generalization ability.
Let $\bm{u}_{\text{test}} = (u_{\text{test}}(y_1^{(\text{test})}),\dots, u_{\text{test}}(y_{M_\text{test}}^{(\text{test})}))^\top$ 
be the discretization of $u_{\text{test}}$,
which is not available in practice.
We measure 
the generalization ability of the DeepONet by means of 
the relative $\ell_2$ error defined by
\begin{equation} \label{def:rel_l2_err}
    \mathcal{E}_{\text{rel}}(\ONet[\bm{f}_{\text{test}}]) := \frac{\sqrt{
    \sum_{i=1}^{M_\text{test}}
    \left(\ONet[\bm{f}_{\text{test}}](y_i^{(\text{test})}) - {u}_{\text{test}}(y_i^{(\text{test})})\right)^2}}{\| \bm{u}_{\text{test}}\|_2}.
\end{equation}

\textbf{Conditional optimality.}
If $\bm{u}_{\text{test}}$ were known, 
by fixing the trunk network,
one can obtain the optimal value $\bm{a}_{\text{test}}^*$
for the branch network at $\bm{f}_{\text{test}}$ by solving
\begin{align*}
    \bm{a}_{\text{test}}^* = \begin{cases}
        \argmin_{\bm{a} \in \mathbb{R}^{N+1}} \| \bm{\Phi}_{\text{test}}(\bm{\mu}^*)\bm{a} - \bm{u}_{\text{test}}\|_2, & \text{with the monolithic method}, \\
        \argmin_{\bm{a} \in \mathbb{R}^{N+1}} \| \bm{\Phi}_{\text{test}}(\bm{\mu}^*)T^*\bm{a} - \bm{u}_{\text{test}}\|_2, & \text{with the two-step method},
    \end{cases}
\end{align*}
where $\bm{\Phi}_{\text{test}}(\bm{\mu}^*)$
is the matrix 
whose $i$-th row is
$\bm{\phi}^\top(y_i^{(\text{test})};\mu^*)$.
Let
\begin{equation} \label{def:cond-opt-rel}
    \ONet^*[\bm{f}_{\text{test}}](y):= \begin{cases}
        \bm{\phi}^\top(y;\mu^*)\bm{a}^*_{\text{test}}, & \text{with the monolithic method}, \\
        \bm{\phi}^\top(y;\mu^*)T^*\bm{a}^*_{\text{test}},   & \text{with the two-step method}.
    \end{cases}
\end{equation}
We then define $\mathcal{E}_{\text{rel}}(\ONet^*[\bm{f}_{\text{test}}])$
as the \textit{optimal} relative $\ell_2$ error.
Here the optimality shall be understood as
conditional 
in the sense that given $\{y_i^{(\text{test})}\}$ and the trunk network $\bm{\phi}(\cdot;\mu^*)$,
$\ONet^*[\bm{f}_{\text{test}}]$
is the least square approximation to $u_\text{test}$.
% which utilizes the optimal branch value.
However, this optimality is not available in practice as $\ONet^*$ requires the target function $\bm{u}_{\text{test}}$ to obtain $\bm{a}^*_{\text{test}}$.
%{\color{red}(Lee: Do we need to define $\bm{a}$?)}

\subsection{Forward Problem: Nonlinear Conductivity}
Let us consider a specific case of 
\eqref{eqn:Darcy}.
Let $f=1$, $g=\cos(x)$, $\partial \Omega = \Gamma_D$
and the conductivity coefficient be $\alpha(p) = \kappa p$ where $\kappa$ is a constant function.
% in the domain $\Omega$. Here, $f = 1$ and we consider only the following Dirichlet boundary condition, $g = \cos(x)$ on $\Gamma$, i.e $\Gamma= \Gamma_D$.
The operator $\mathcal{G}$ of interest is
\begin{equation*}
    \mathcal{G}: \mathcal{X} \ni \kappa(\cdot) \mapsto p(\cdot) \in \mathcal{Y},
\end{equation*}
where $\mathcal{X} =\{\kappa \ | \  \kappa(x,y) = \beta, \forall (x,y) \in \Omega,  \beta \in [1,1000] \}$
and
$\mathcal{Y}$ is an appropriate space 
where the solution $p$ lies.
% {\color{red} YJ: Can you name a specific space of $\mathcal{Y}$? E.g. $L_2(\Omega)$?}
Note that for any $\kappa \in \mathcal{X}$,
it is well-known \cite{evans2022partial} that there exists a unique solution $p(\cdot)$ of the system \eqref{eqn:Darcy}.  
% In this case, $\Omega_{x} = \Omega_{y} = \Omega$.
% e are interested in approximating the solution operator $\mathcal{G}$. 
  
Since the input functions are constant functions,
we simply set $\mx = 1$, e.g., $x_1=(0,0)$.
% Thus, the discretized input dimension is oe $(\mx =1)$ for this case.
Accordingly, the input data are generated as the collection of 1000 equidistant $\beta$ values in $[1,1000]$, i.e., 
$\{1,2,\dots,1000\}$.
The corresponding output data are obtained by the FEM solver on 2049 grid points, i.e., $\my=2049$.
The data is then randomly split into two -- 
900 of them are used as training and the remaining 100 are used as test data. 
We employ a DeepONet whose branch and trunk architectures are
$\vec{\bm{n}}_b = (1,500,51)$
and $\vec{\bm{n}}_t = (2, 50, 50, 50, 50)$, respectively. Both branch and trunk networks use the rectified linear unit (ReLU) activation function and were initialized by the He initialization scheme \cite{he2015delving}.
Throughout, we employ \texttt{Adam} optimizer \cite{kingma2014adam} with full-batch.

\begin{figure}[htbp]
	\centerline{
	    \subfloat[ ]{\label{fig:ex1_loss:a}
	    \includegraphics[width=6.3cm]{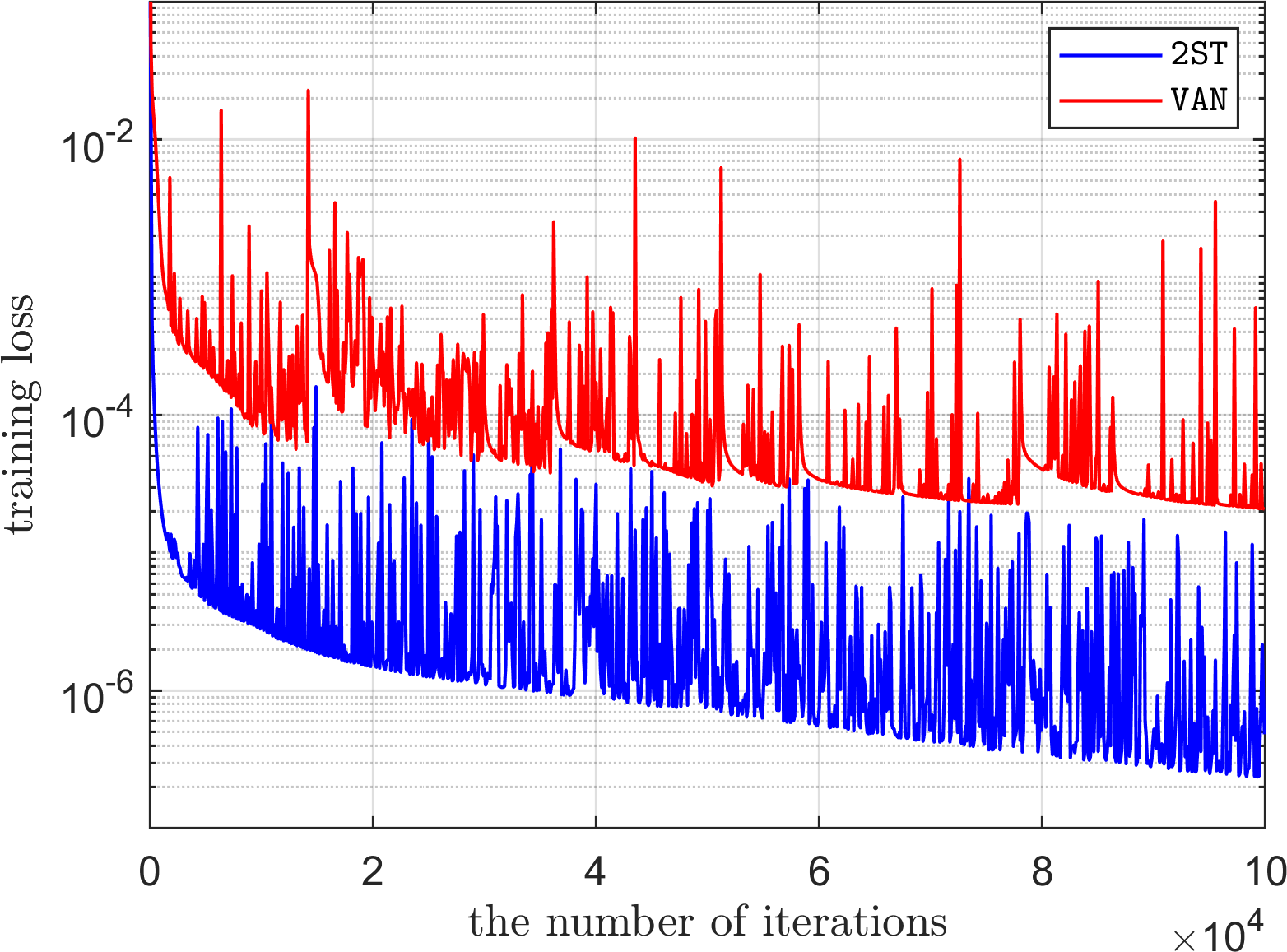}
	    }
	    \subfloat[ ]{\label{fig:ex1_loss:b}
		\includegraphics[width=6.3cm]{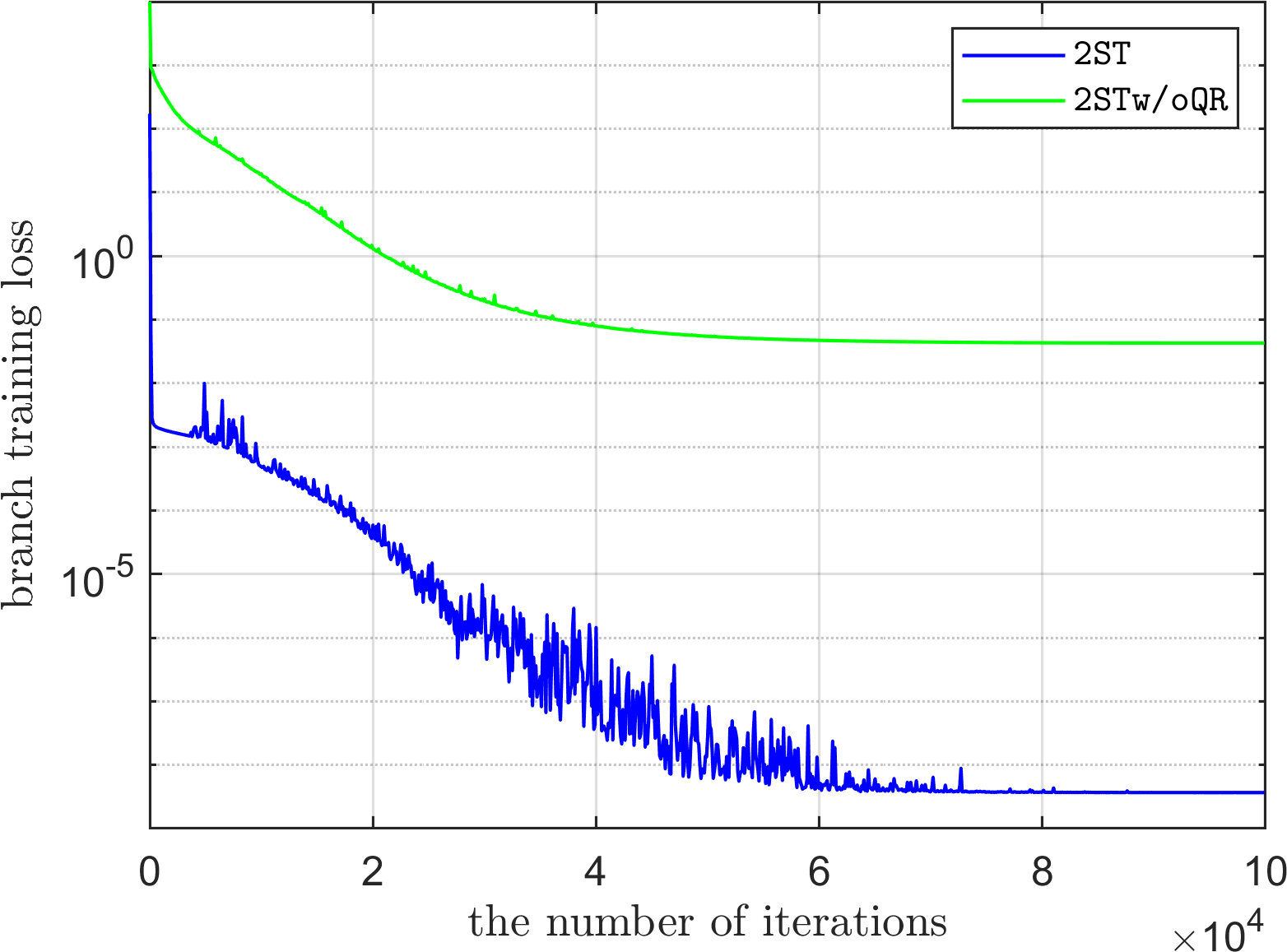}
	    }}
	\caption{Example 5.1.
		% The training loss versus the number of iterations by \texttt{2ST}, \texttt{VAN}, \texttt{2STw/oQR}.
        (Left)
        The training loss versus the number of iterations. Here the loss refers to \eqref{def:loss} for \texttt{VAN}
        and \eqref{train-trunk}
        for \texttt{2ST} (thus \texttt{2STw/oQR}). 
        (Right)
        The branch loss \eqref{train-branch} are reported for \texttt{2ST} and \texttt{2STw/oQR}. This shows the effectiveness of orthogonalization in the second step of the proposed training method.
	}
	\label{fig:ex1_loss}
\end{figure}

In Figure~\ref{fig:ex1_loss:a}, we plot the training loss 
versus the number of iterations
by both \texttt{2ST} and \texttt{VAN}.
Specifically, the training loss refers to 
the trunk network loss which is defined in \eqref{train-trunk} for \texttt{2ST}  
and the standard overall loss of \eqref{eqn:DeepONet-train}
for \texttt{VAN}.
It can be clearly seen that 
the loss by \texttt{2ST} is roughly two orders of magnitude smaller than the one by \texttt{VAN}.
As a matter of fact, the smallest loss attained by \texttt{2ST} is $2.33 \times 10^{-7}$,
while the one by \texttt{VAN} is $2.08\times 10^{-5}$.
This is not a single isolated case.
We tested five independent simulations, and the averaged smallest loss achieved by the two methods is $2.11 \times 10^{-7}$ and $1.84 \times 10^{-5}$ for \texttt{2ST} and \texttt{VAN} respectively.
This demonstrates the effectiveness of the proposed two-step method for learning the trunk network.
The remaining task for \texttt{2ST} is then to learn the branch network following \eqref{train-branch}.

In Figure~\ref{fig:ex1_loss:b}, the training loss for the branch network versus the number of iterations is plotted.
We can see that the loss reaches the level of $10^{-9}$ by \texttt{2ST}
which utilizes QR-factorization,
while the one without QR cannot reach a similar level.
This demonstrates the effectiveness of the orthogonalization in the proposed two-step training method.

\begin{figure}[!h]
\centering
{\includegraphics[width=0.48\textwidth]{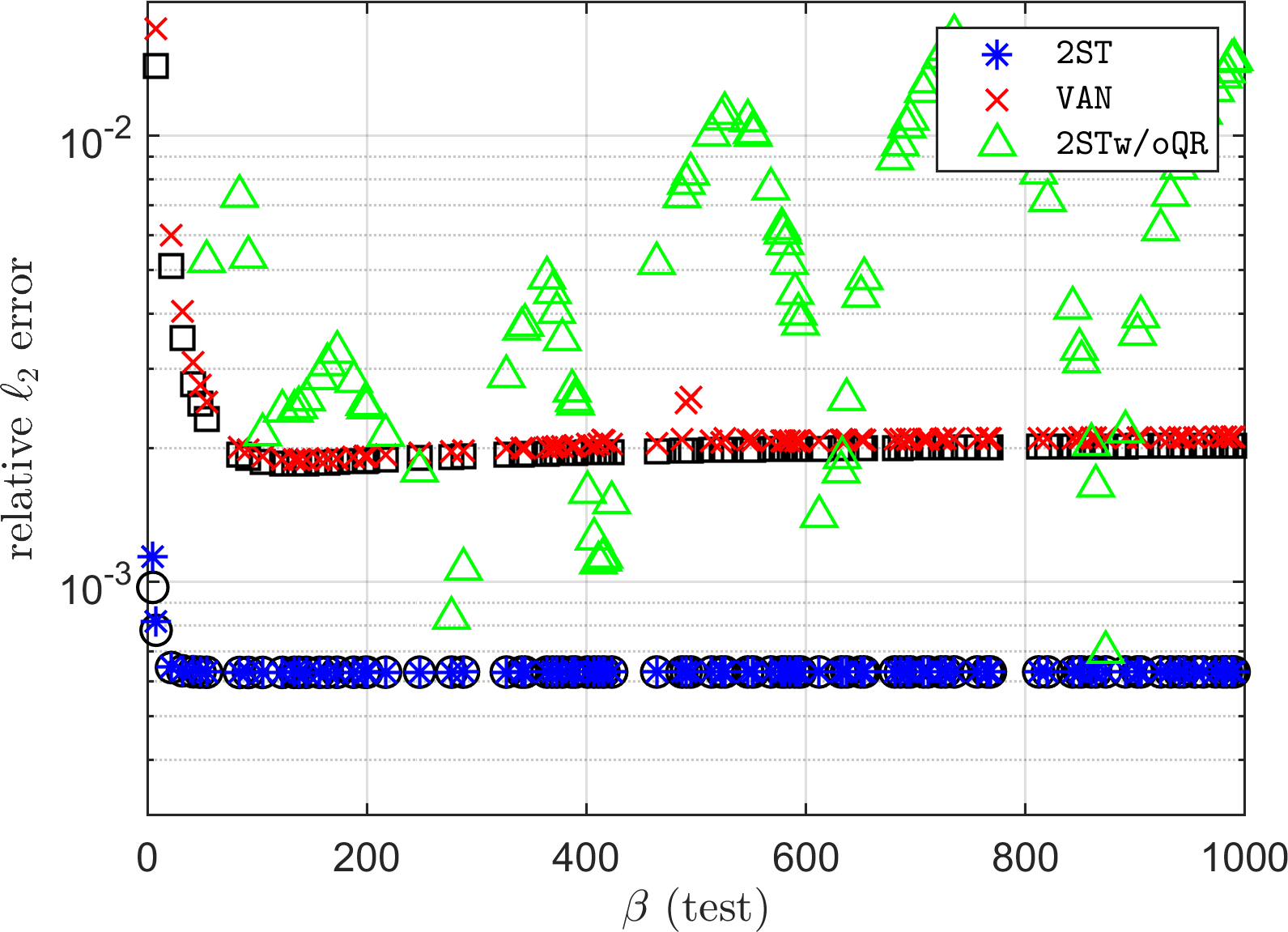}
\includegraphics[width=0.48\textwidth]{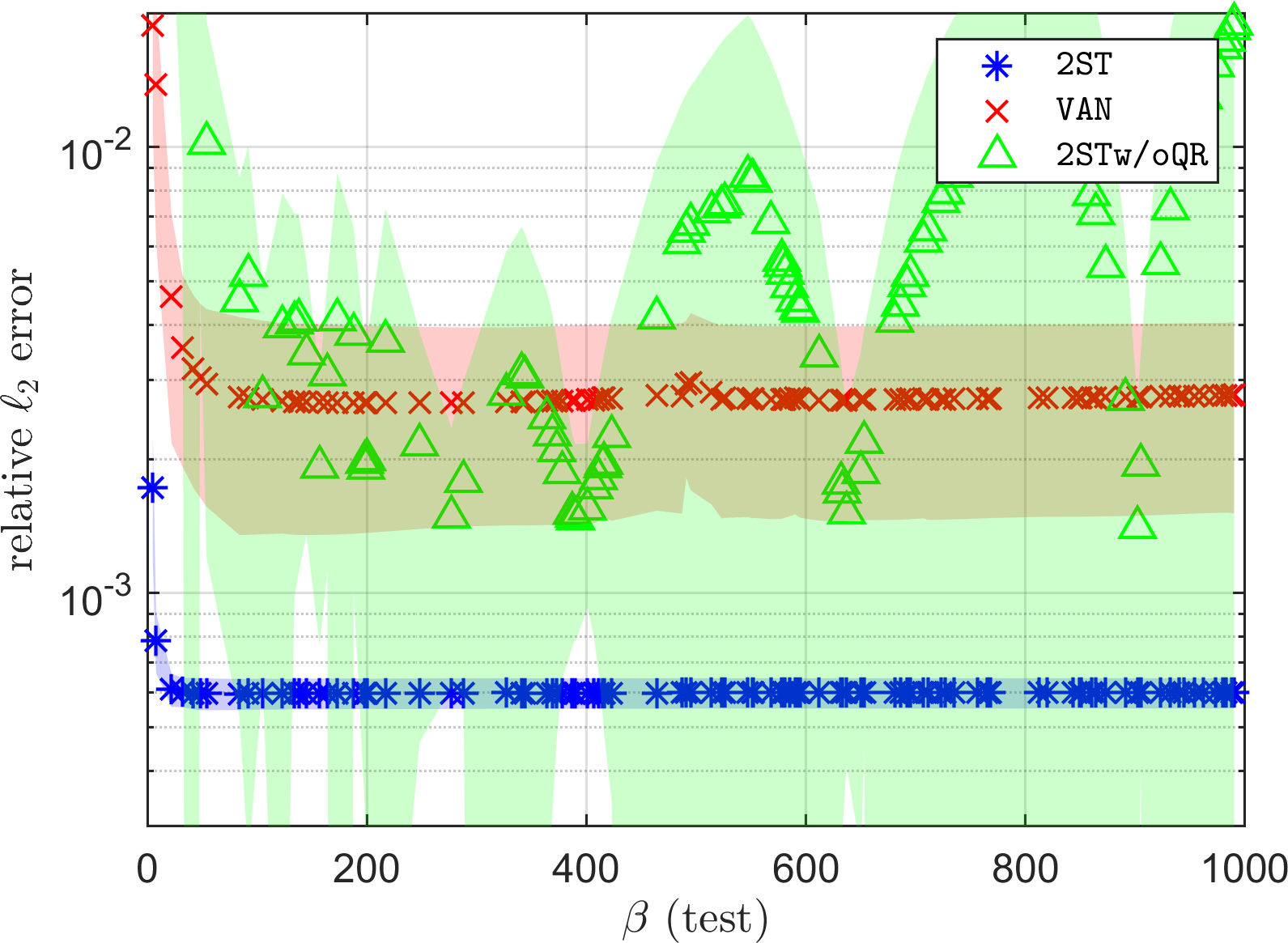}} \\[10pt]
{\includegraphics[width=0.325\textwidth]{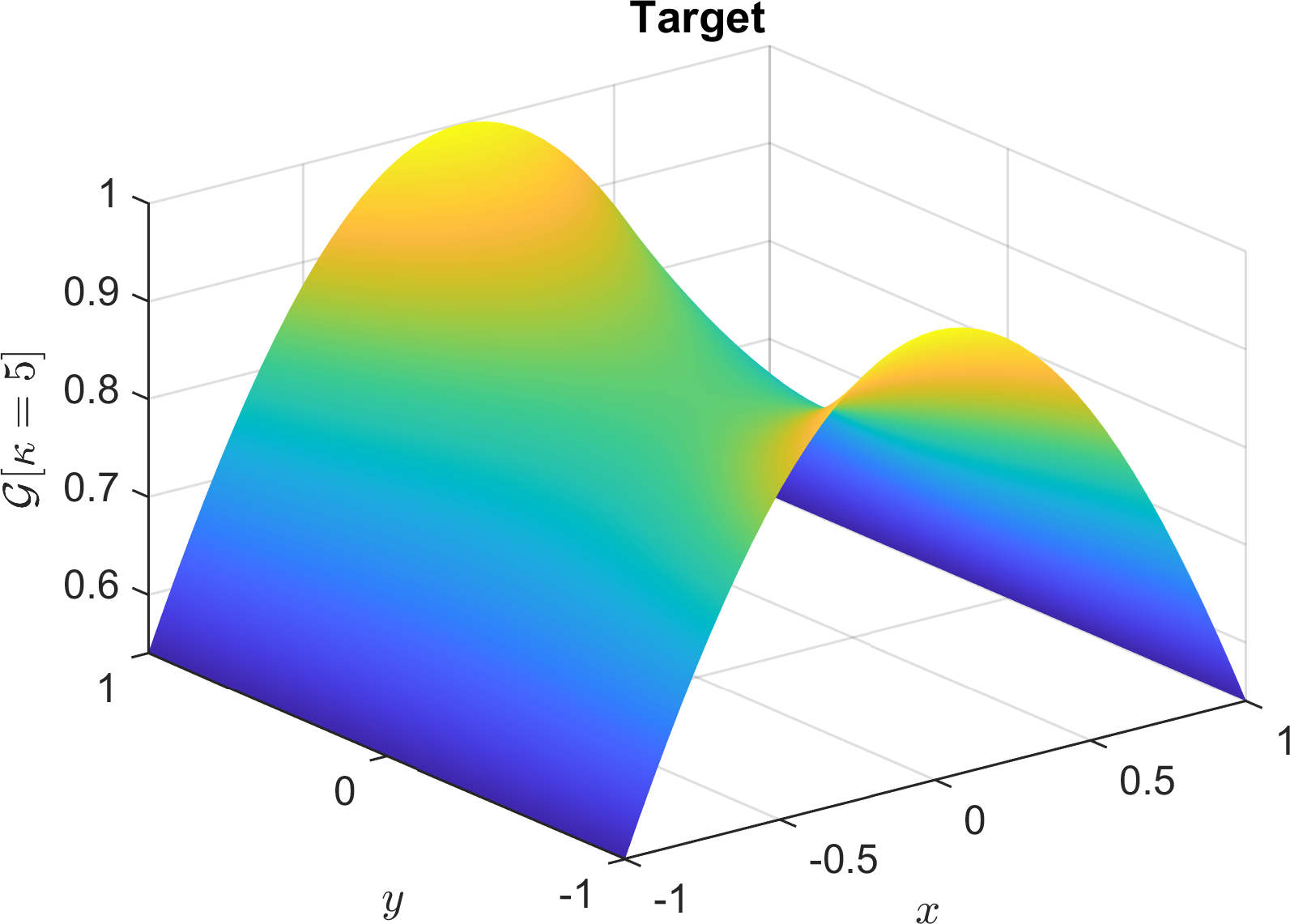}
\includegraphics[width=0.325\textwidth]{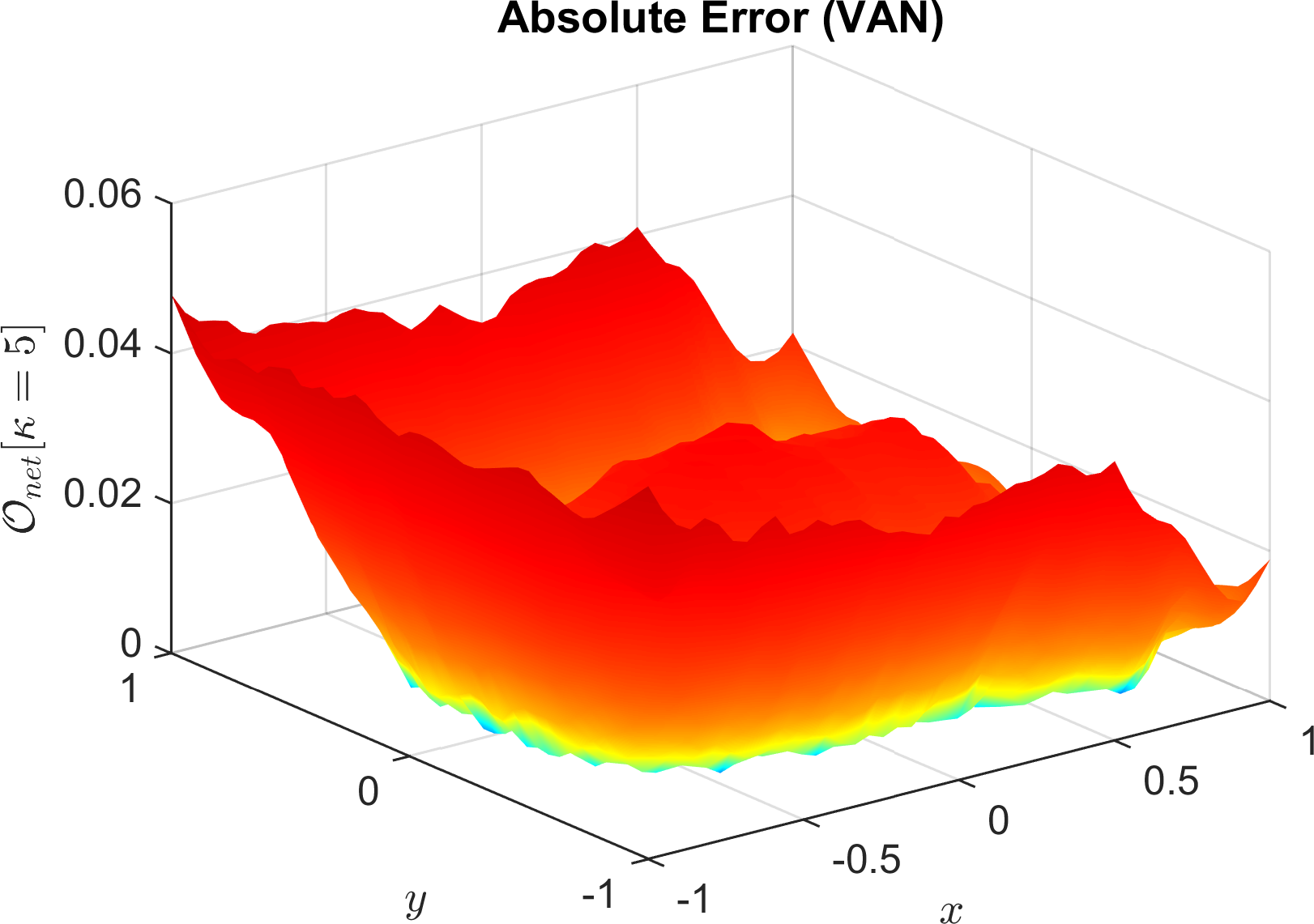}
\includegraphics[width=0.325\textwidth]{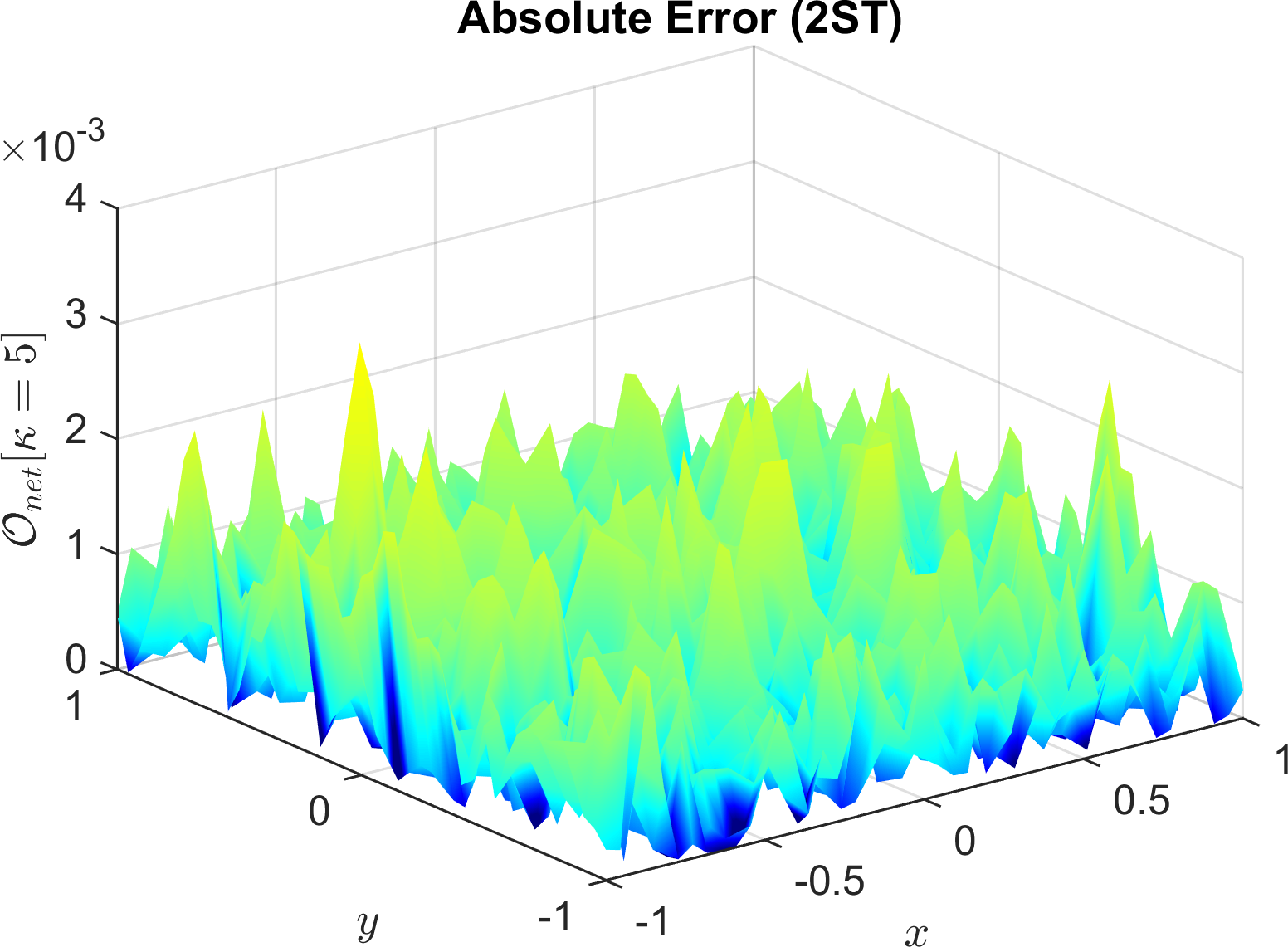}
}
\caption{Example 5.1.  (Top) The relative $\ell_2$ errors versus the 100 test $\beta$ values
by \texttt{2ST}, \texttt{VAN} and \texttt{2STw/oQR}.
(Bottom) Left: The graph of the target operator evaluated at a test function $\kappa = 5$,
which is the smallest value lying in the test data. 
Middle and Right: 
The absolute error maps by \texttt{VAN} (middle) and \texttt{2ST} (right).
}
\label{fig:ex1_test}
\end{figure}

In the top left of Figure~\ref{fig:ex1_test},
the relative $\ell_2$ errors defined in \eqref{def:rel_l2_err} are
reported. 
The errors are plotted with respect to the 100 test $\beta$ values, where the $\beta$ are constant input functions.
The results of \texttt{2ST}, \texttt{VAN}, \texttt{2STw/oQR}
are marked as 
asterisks ($\ast$),
crosses ($\times$),
triangle ($\triangle$),
respectively. 
The optimal relative $\ell_2$ errors defined in \eqref{def:cond-opt-rel} are shown as circles ($\circ$) for \texttt{2ST} and squares ($\square$) for \texttt{VAN}, and they serve as references.
%%%%%
We note that having small optimal relative errors means the trunk network was well-trained to represent unseen output functions. In a similar vein, it is also worth noting that if the relative errors of \texttt{VAN} (\texttt{2ST}) are close to the optimal ones of \texttt{VAN} (\texttt{2ST}), it implies that the branch network was well-trained, allowing the DeepONet to generalize effectively.
%%%%%
It is clearly observed that 
\texttt{2ST} achieves the smallest relative $\ell_2$ errors being almost identical to the optimal ones, while \texttt{2STw/oQR} yields much higher and unstable errors.
This again demonstrates the effectiveness of orthogonalization in the proposed method.
\texttt{VAN} is able to produce relative errors similar to optimal ones,
however, the relative optimal errors of \texttt{VAN} are much higher than the ones by \texttt{2ST}.
This indicates the ineffectiveness of the monolithic training for learning the trunk network.
On contrary,
the proposed two-step training method \texttt{2ST} effectively trains 
not only the trunk network but also the branch network.
On the top right, we report the means of the relative $\ell_2$ errors from five independent simulations.
The shaded area is the area that falls within one standard deviation of the mean.
On the bottom left, we report the graph of the target output function for $\beta=5$ (the smallest value belongs to the test data).
On the bottom middle and right,
the absolute error maps by \texttt{VAN} and \texttt{2ST} are shown, respectively.
It is clearly observed that 
the absolute error by \texttt{2ST}
is at least one order magnitude smaller than
the one by \texttt{VAN}.

% \clearpage
% \newpage
\subsection{Inverse Problem: Discontinuous Conductivity}
\begin{figure}[!h]
\centering
\begin{subfigure}[b]{0.33\textwidth}
\centering
\begin{tikzpicture}[scale=1.85]
\filldraw[color=red!60, fill=red!5, very thick](0,0) circle (0.5);
\draw[thick] (-1,-1) rectangle (1,1);

\draw (0.,0.85) node[] {$\Gamma_{D_1}$};
\draw (0.,-1.) node[above=0pt] {$\Gamma_{N_1}$};
\draw (-1,0.) node[right=0pt] {$\Gamma_{N_2}$};
\draw (0.8,0.) node[] {$\Gamma_{N_2}$};

   \draw (-1,1)  node[below=10pt, right=5pt] {$\Omega_2$};
   \draw (0,0) node[above=10pt] {$\Omega_1$};

   \draw (-1,-1)  node[above=20pt,right=5pt] {$\kappa_2 = 1$};
   \draw (-0.45,-0.6) node[above=20pt,right=5pt] {$\kappa_1 = \beta$};
   
\end{tikzpicture}
\caption{Domain Geometry}
\label{fig:domain}
\end{subfigure}
% \hfill
\begin{subfigure}[b]{0.3\textwidth}
         \centering
         \includegraphics[width=\textwidth]{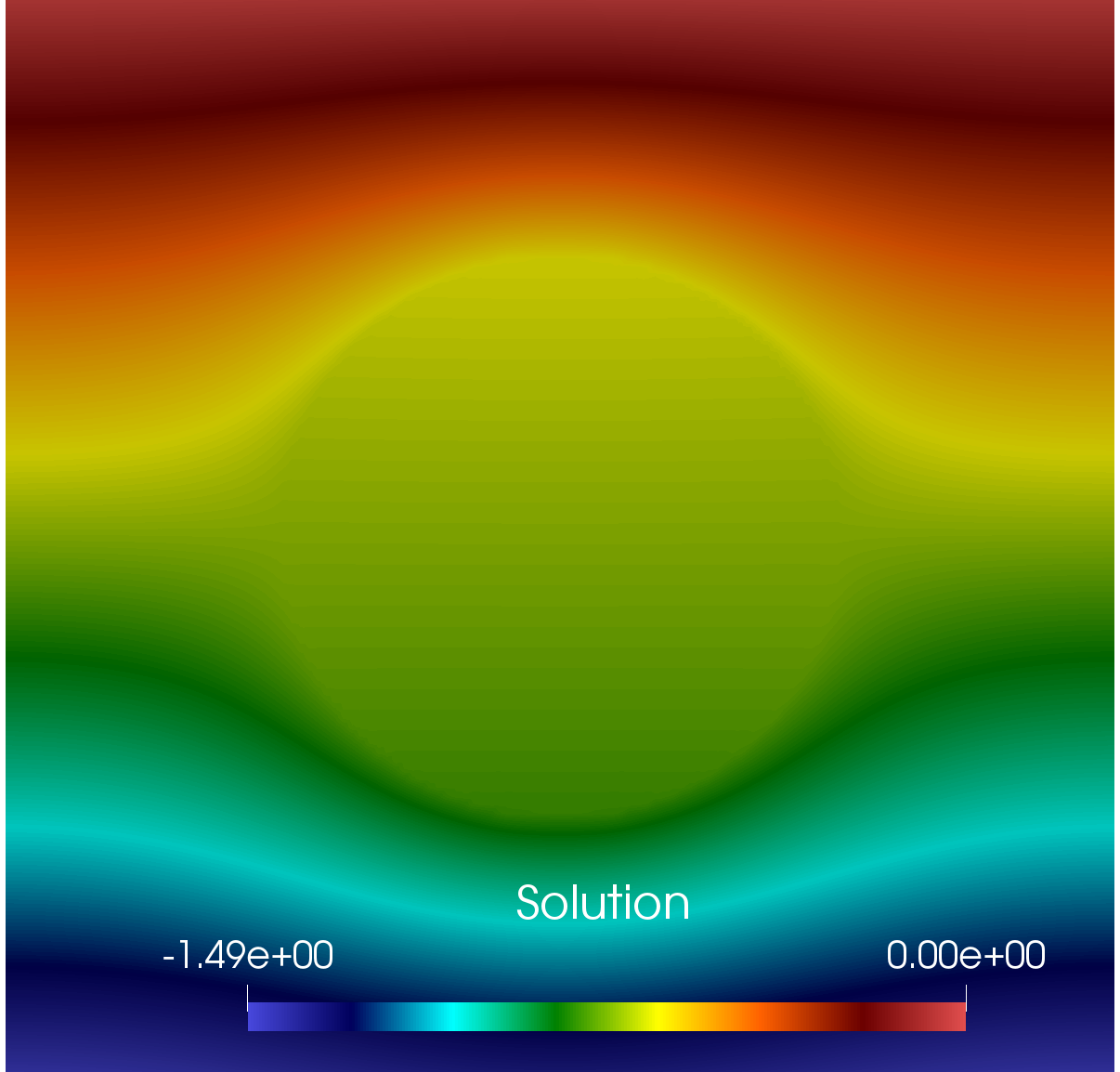}
         \caption{$p(\cdot)$}
         \label{fig:k1}
\end{subfigure}
\begin{subfigure}[b]{0.35\textwidth}
         \centering
         \includegraphics[width=\textwidth]{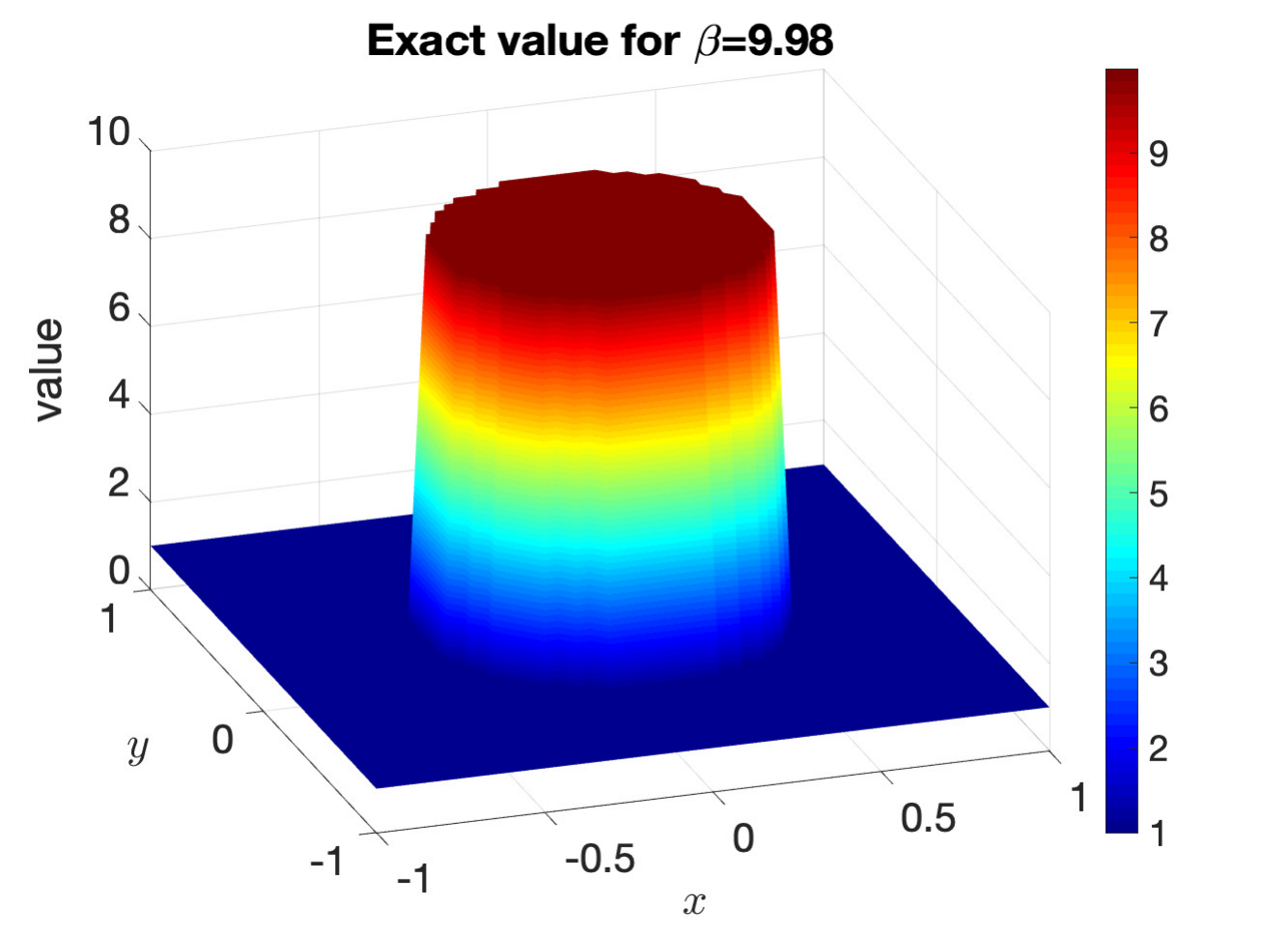}
         \caption{$\kappa(\cdot;\beta= 9.98)$}
         \label{fig:k2}
\end{subfigure}
% \hfill
\caption{Example 2. (a) Problem domain with detailed geometries. 
(b) The solution $p$ of \eqref{eqn:Darcy} with the boundary conditions \eqref{eqn:2_bc} at $\kappa(\cdot;\beta=9.98)$.
(c) The graph of $\kappa(\cdot;\beta=9.98)$.}
\label{fig:ex_1_2}
\end{figure}
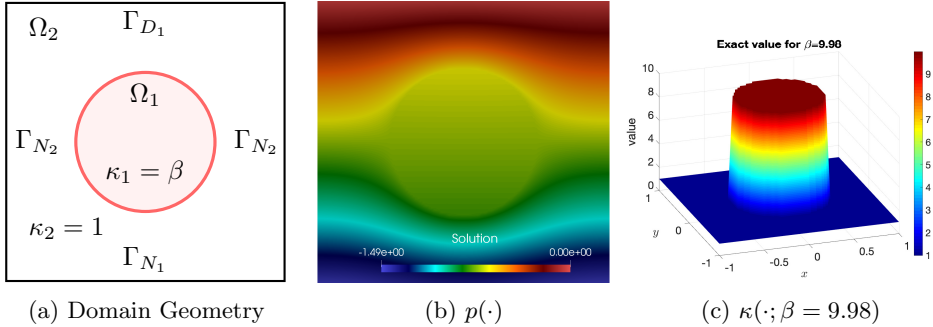

In this case, let $f=0$ and consider a piece-wise constant (discontinuous) conductivity  
$\alpha := \kappa(\cdot;\beta)$
where $\kappa$ is defined by 
\begin{equation}
\kappa(x,y; \beta)= 
\begin{cases}
\beta  & \text{ if } (x, y) \in \Omega_1, \\
1     & \text{ if }  (x, y) \in \Omega \backslash \Omega_1,
\end{cases}
\label{eqn:2_kappa}
\end{equation}
%and test the algorithm with discontinuities in the input value $\kappa := \kappa(x,y;\beta)$. 
% In particular, the conductivity coefficient 
%\begin{equation} \label{ex1.1_beta}
%\kappa_1 := \kappa|_{\Omega_1} = \beta \text{ in } \Omega_1, \ \    
%\kappa_2 := \kappa|_{\Omega_2} = 1 \text{ in } \Omega_2. 
%\end{equation}
where $\Omega_1$ is a disk centered at the origin $(0,0)$ with radius $r=0.5$.
% and  $\Omega_2 = \Omega \backslash \Omega_1$.
%Thus, we test the algorithm with discontinuities in the input value $\kappa := \kappa(x,y;\beta)$. 
The Dirichlet and Neumann boundary conditions are imposed by
\begin{equation}
    p = 0 \text{ on } \Gamma_{D_1}, \ \ 
    -\kappa \nabla p \cdot {\bf{n}} = 0 \text{ on } \Gamma_{N_2}, \ \ 
    -\kappa \nabla p \cdot {\bf{n}} =  1 \text{ on } \Gamma_{N_1}.
\label{eqn:2_bc}    
\end{equation}
Figure \ref{fig:domain} shows the detailed geometry of the problem.

We are interested in the inverse problem of \eqref{eqn:Darcy} with the boundary conditions of \eqref{eqn:2_bc}.
That is, the operator $\mathcal{G}$ of interest maps 
a given solution $p$
to the corresponding conductivity coefficient $\kappa$ as a function.
Specifically, 
$\mathcal{G}:\mathcal{X} \ni p(\cdot) \mapsto \kappa(\cdot) \in \mathcal{Y}$
where 
$\mathcal{X}$ is the collection of solutions 
obtained at various $\kappa(\cdot;\beta)$ 
defined in \eqref{eqn:2_kappa}
where $\beta \in [0.01,10]$ for $\kappa$,
and $\mathcal{Y}$ is the collection of the corresponding $\kappa(\cdot;\beta)$. 
Figures \ref{fig:k1} and~\ref{fig:k2}
show an input-output pair for the operator.
% Figure \ref{fig:k2}  plots the $\kappa(x,y;\beta)$ where $\beta=9.98$, and 
% Figure \ref{fig:k1} illustrates the corresponding  approximated solution of $p(x,y)$.
%a different choice for trunk net training. 
%\begin{wrapfigure}{r}{0.35\textwidth}
%\end{wrapfigure}

For data generation, 
we consider the collection of 1000 $\beta$ values -- $\{\beta = j \times 0.01: j \in \{1,2,\dots, 1000\}\}$ -- and 
solve the corresponding equation \eqref{eqn:Darcy} to obtain $p$.
This is done by the FEM solver on $\mx = 4225$ points in $\Omega$. 
On that exact grid, i.e., $\my = 4225$, the output function data are generated
according to \eqref{eqn:2_kappa}.
The data are then randomly split into 900 training data and 100 test data.
We employ a DeepONet with the trunk and branch architectures of 
$\vec{\bm{n}}_t=(2, 25, 25, 25, 25, 25)$
and 
$\vec{\bm{n}}_b=(4225, 100, 26)$, respectively,
with the ReLU activation function.
Both are initialized with the He scheme \cite{he2015delving}.
Since the effectiveness of QR was already demonstrated in the previous example, here we only consider $\texttt{2ST}$.

\begin{figure}[htbp]
	\centerline{
	    \subfloat[ ]{\label{fig:ex2_loss:a}
	    \includegraphics[width=6.3cm]{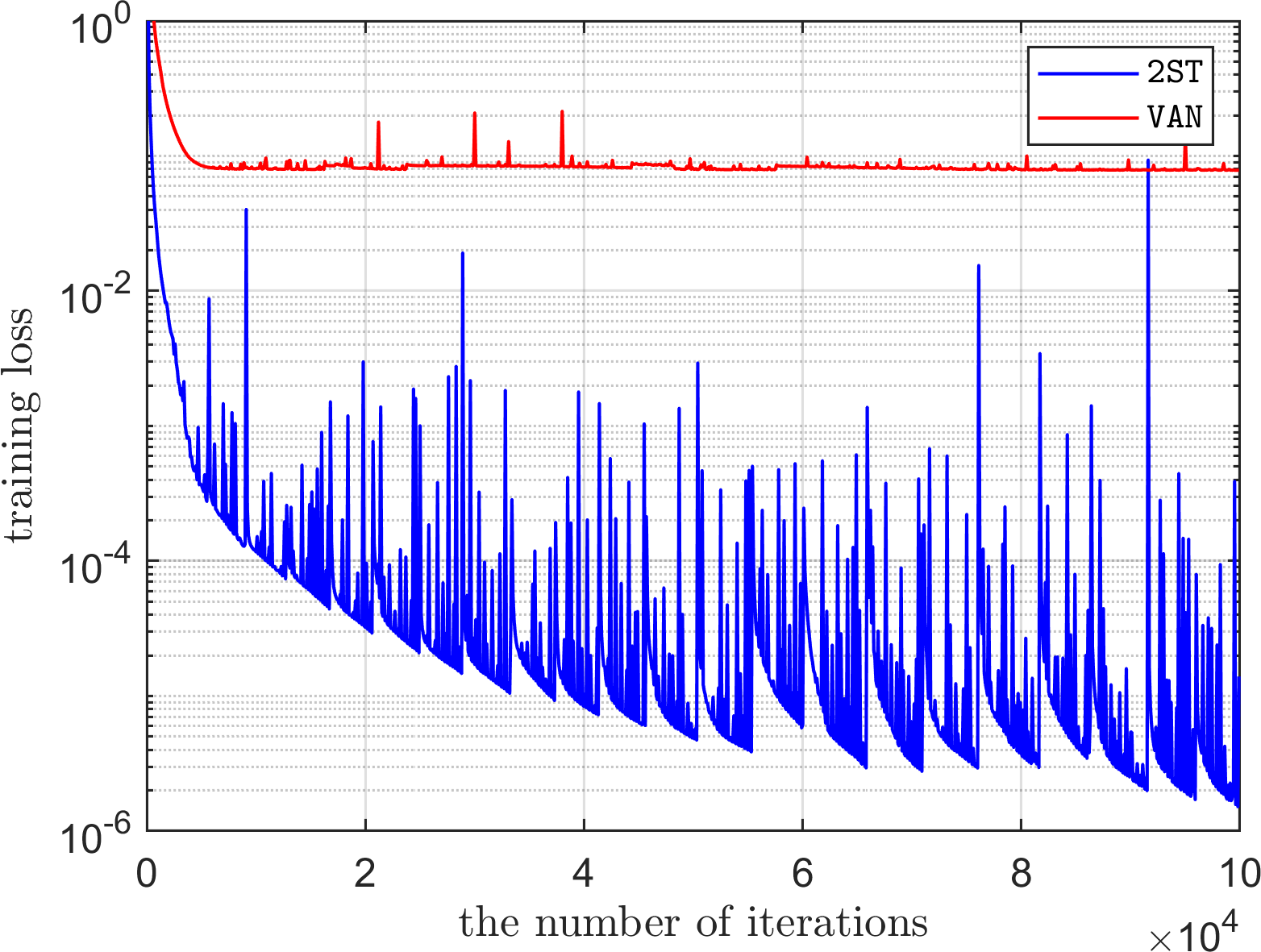}
	    }
	    \subfloat[ ]{\label{fig:ex2_loss:b}
		\includegraphics[width=6.3cm]{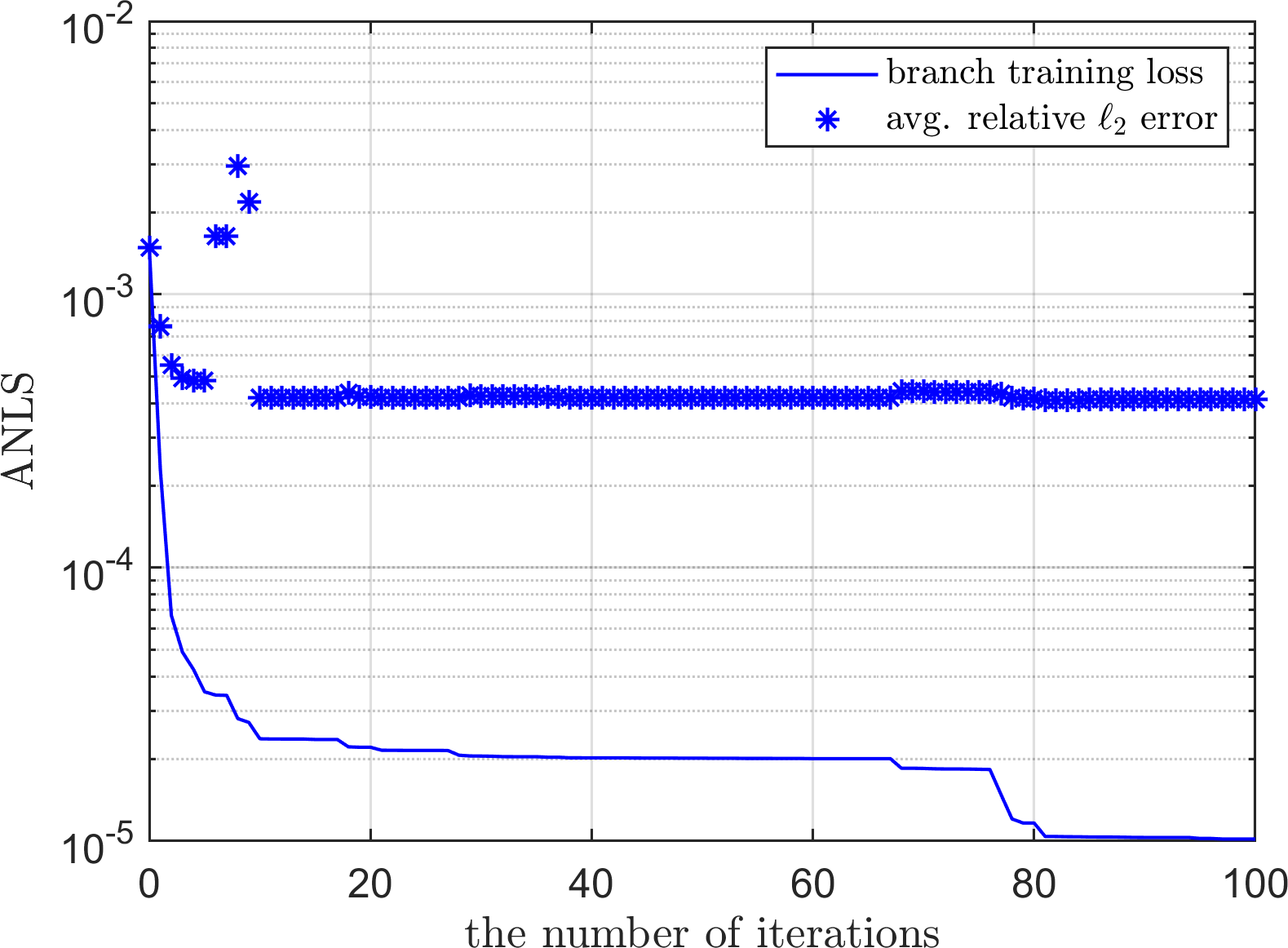}
	    }}
	\caption{Example 5.2.
		% The training loss versus the number of iterations by \texttt{2ST}, \texttt{VAN}, \texttt{2STw/oQR}.
        (Left)
        The training loss versus the number of iterations. Here the loss refers to \eqref{def:loss} for \texttt{VAN}
        and \eqref{train-trunk}
        for \texttt{2ST}. 
        (Right)
        The branch training loss of \texttt{2ST} versus 
        the number of ANLS \cite{ainsworth2021plateau,ainsworth2022ANLS} iterations. 
        Also, the average relative $\ell_2$ errors over the 100 test data are marked as asterisks ($\ast$).
	}
	\label{fig:ex2_loss}
\end{figure}

Figure~\ref{fig:ex2_loss:a} shows the training loss versus the number of iterations by both \texttt{2ST} and \texttt{VAN}. Specifically, it displays the trunk network training loss for \texttt{2ST} and the overall standard loss for \texttt{VAN}.
We employ \texttt{Adam}  optimizer \cite{kingma2014adam} with full-batch and its default hyperparameters.
It is clearly observed that
the two-step training method \texttt{2ST}
minimizes the loss to the level of $10^{-6}$,
while the monolithic standard training \texttt{VAN} stagnates at the level of $10^{-1}$ after $10,000$ iterations.
This again indicates that \texttt{2ST} effectively trains the trunk network
to represent the output functions, which are the piece-wise constant functions \eqref{eqn:2_kappa}.

For the second step \eqref{train-branch} of \texttt{2ST}, we employ the Active Neuron Least Squares (ANLS) training method 
developed in \cite{ainsworth2021plateau,ainsworth2022ANLS}.
This is possible because the second step is merely a standard regression task on which ANLS is applicable.
In Figure~\ref{fig:ex2_loss:b},
the branch training loss is plotted
with respect to the number of ANLS iterations.
ANLS minimizes the branch loss to the level of $1.02 \times 10^{-5}$ merely within 100 iterations.
The averaged relative $\ell_2$ errors over the 100 test data are also reported.
It can be seen that the average test errors are saturated after merely 20 ANLS iterations.
This indicates that the branch network is successfully trained by ANLS.
We remark that ANLS is not applicable to the monolithic training of DeepONets.
% We remark that other optimization methods (e.g. \texttt{Adam} \cite{kingma2014adam}, \texttt{L-BFGS}) are still applicable
% if the architecture of the branch networks is sufficiently large enough.

% {\color{red}Figure \ref{fig:ex2_2} illustrates the relative $l_2$ errors for each test cases for $\kappa$. The comparison between monolithic \texttt{Adam}, two-step \texttt{Adam}, and two-step ANLS emphasizes the importance of the training algorithm, and we observe the two-step algorithms perform better the monolithic training, and the two-step ANLS shows the best performance. 
% }

% {\color{red}Finally, Figure \ref{fig:ex2_2_trunk} a)-c) presents some examples of trunk networks, $\phi_j(\cdot, \mu_j)$ where $j=1,2,$ and $3$. 
% We observe that these trunk net values are constructed by detecting the given discontinuity. }

% The challenge in this example is to consider the discontinuity in the coefficient $\kappa$ and we test the capability of our proposed algorithm with 
% the discontinuities in the input value.

\begin{figure}[htbp]
	\centering
	\includegraphics[width=0.7\textwidth]{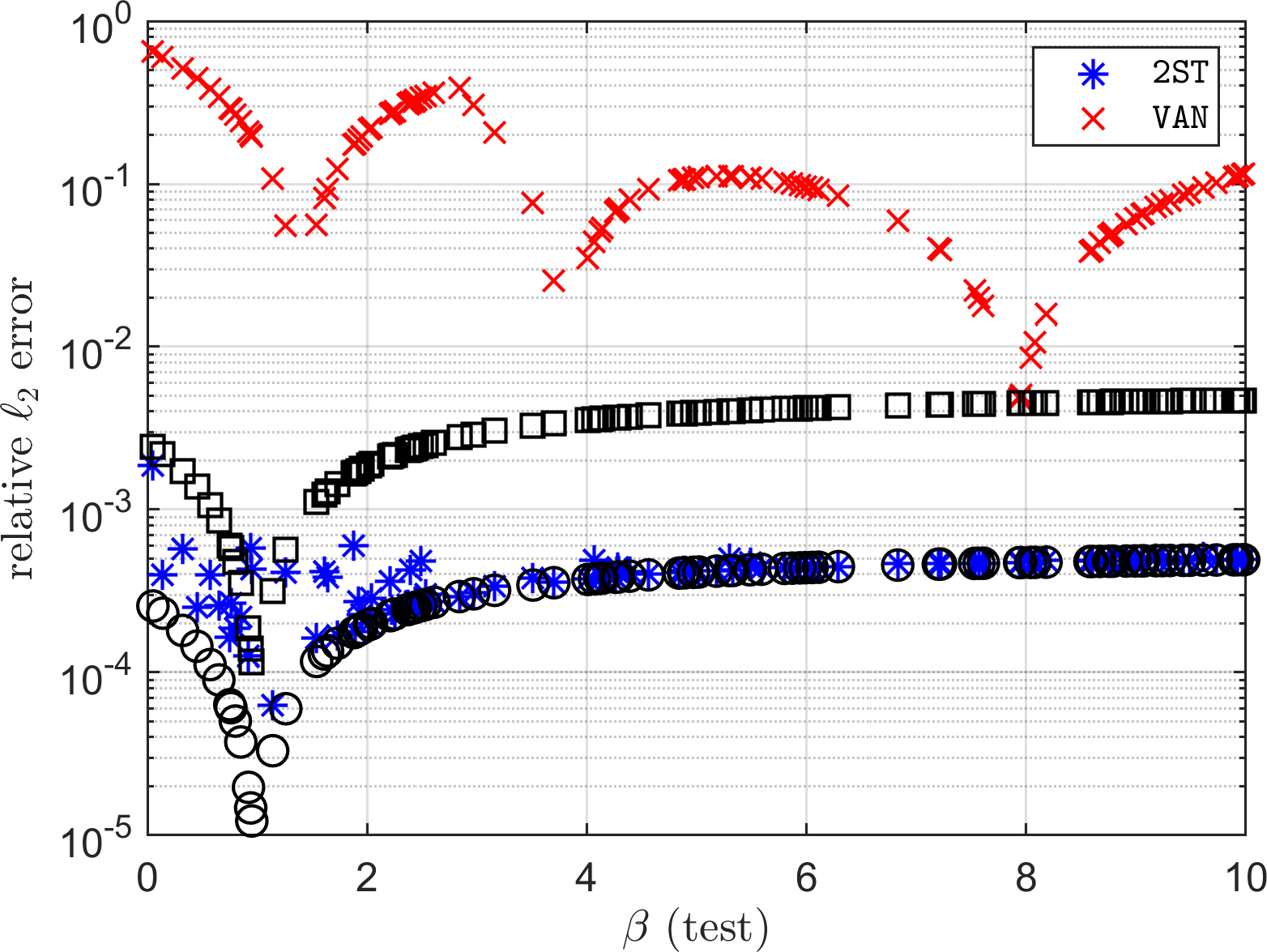}
	     \\[10pt]
	\includegraphics[width=0.49\textwidth]{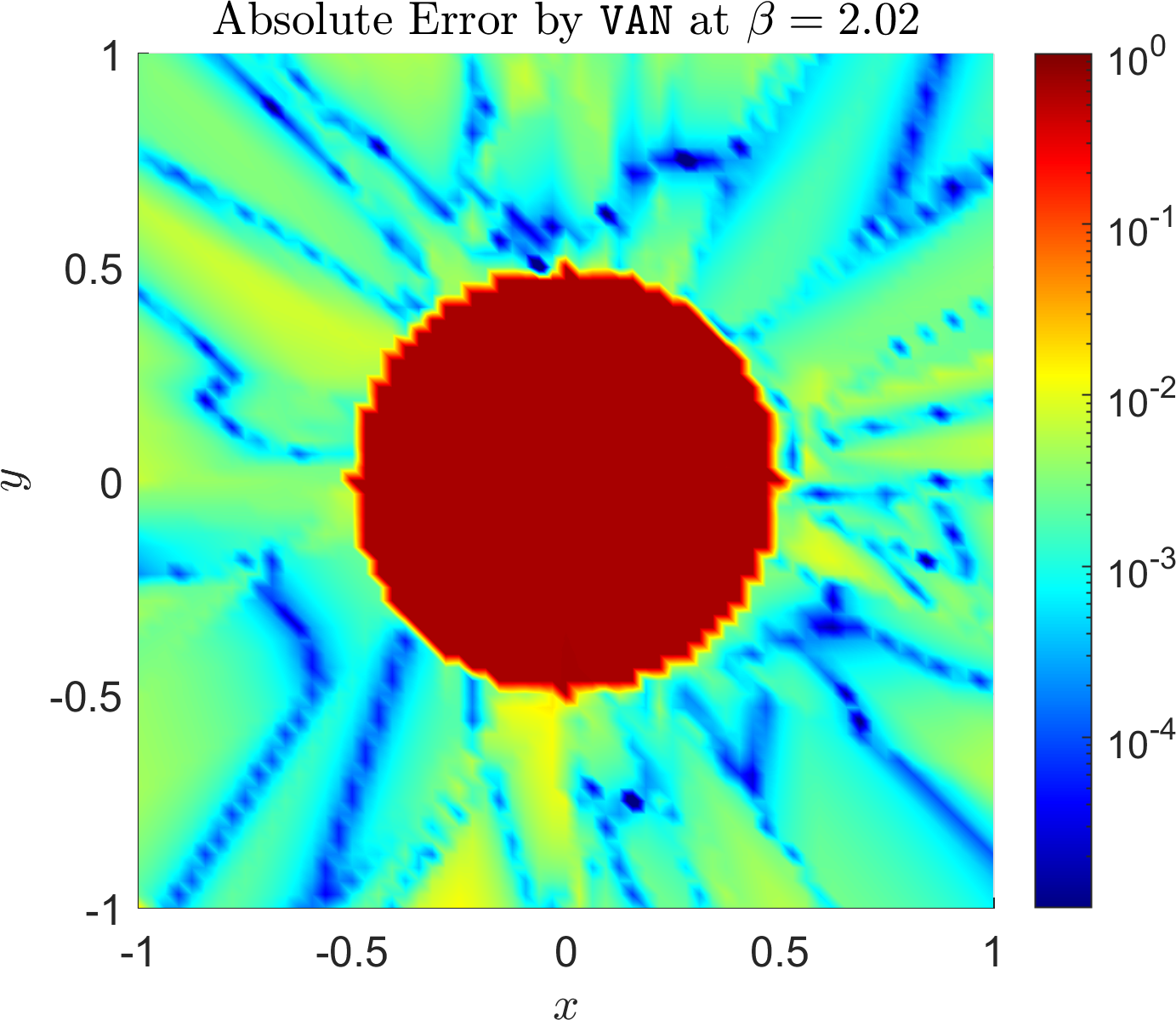}
        \includegraphics[width=0.49\textwidth]{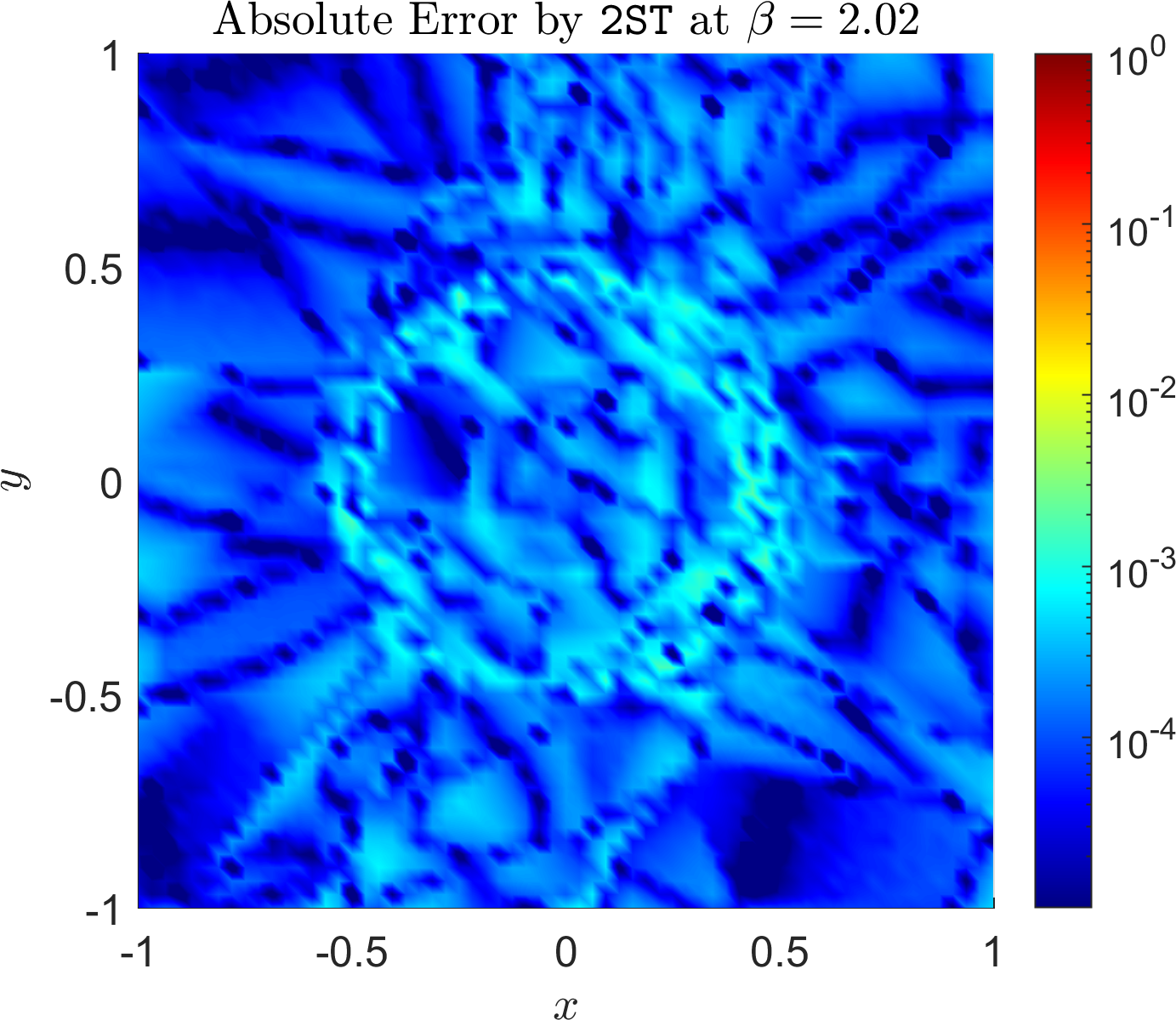}
	\caption{Example 5.2.
        (Top) The relative $\ell_2$ errors on the test data are reported for \texttt{2ST} and \texttt{VAN}.
        (Bottom)
        The absolute error maps 
        for $\beta = 2.02$ (test)
        by \texttt{VAN} (left)
        and \texttt{2ST} (right).
	}
	\label{fig:ex2_test}
\end{figure}

At the top of Figure~\ref{fig:ex2_test}, we plot the relative $\ell_2$ errors for the 100 test data. Since each test datum is determined by its corresponding $\beta$ value, the errors are plotted with respect to the test $\beta$ values. The optimal relative $\ell_2$ errors \eqref{def:cond-opt-rel} are also reported and serve as reference benchmarks. The optimal errors for \texttt{VAN} and \texttt{2ST} are indicated by squares ($\square$) and circles ($\circ$), respectively.
Again note that these optimal values are not available in practice as the underlying target output function data are required. 
Therefore, the closer to the optimal values from the trained DeepONet, 
the better generalization performance it implies.
It can be seen that the optimal errors by \texttt{2ST} are roughly one order of magnitude smaller than those by \texttt{VAN}.
This indicates that 
\texttt{2ST} can train 
the trunk network more effectively 
than \texttt{VAN}.
Furthermore, it is clearly observed that 
the relative test errors by \texttt{2ST} 
are close to the optimal ones 
and are even almost identical especially if $\beta \in [3,10]$.
On the other hand, 
the test errors by $\texttt{VAN}$ are 
way off from the corresponding optimal ones. 
The averaged relative error over the 100 test data by \texttt{VAN} and \texttt{2ST} are 
$1.48\times 10^{-1}$ 
and 
$4.14\times 10^{-4}$, respectively.
On the bottom of Figure~\ref{fig:ex2_test}, 
the absolute error maps of DeepONets trained by both \texttt{VAN} and \texttt{2ST} at $\beta=2.02$ (test) are shown.
It can be seen that while both capture well the discontinuity, 
\texttt{2ST} can accurately predict the value on the circle. 
This again demonstrates the effectiveness of the proposed two-step training method over the standard monolithic one.

\subsection{Multiple Inputs and Nonlinear Conductivity}
Let us consider the nonlinear conductivity
$\alpha(p):= \kappa p$
and the Dirichlet boundary condition (i.e $\Gamma_N = \emptyset$)
in \eqref{eqn:Darcy}.
In this case, we illustrate the capabilities of the proposed algorithm by considering the solution operator whose inputs are the triplet of 
the right-hand side source term $f$, 
the conductivity $\kappa$, and
the boundary value $g$.
That is, the operator of interest is 
\begin{equation*}
    \mathcal{G}:\mathcal{X} \ni (f,\kappa,g)   \mapsto p \in \mathcal{Y}.
\end{equation*}
The input function space is 
the set of triplets
$\mathcal{X} = \{(f,\kappa,g) | f \in \mathcal{F}, \kappa \in \mathcal{K}, g \in \mathcal{G}\}$
where $\mathcal{F}$, $\mathcal{K}$,
$\mathcal{G}$ are all the collection of constant functions in the range of $[0.1, 10]$.
The output space $\mathcal{Y}$ is 
the collection of the corresponding solutions $p$ to the system \eqref{eqn:Darcy}.
Note that 
for any $(f,\kappa,g) \in \mathcal{X}$,
there exists a unique solution $p$ to \eqref{eqn:Darcy}. 

% In this case, we set the right-hand side source term ($f$), the conductivity $\kappa$, and the boundary values $g$ to be each constant in the range of   $[0.1, 10]$.  
To generate a data set, 
we consider the grid of 1M points in $[0.1,10]^3$,  
\begin{equation*}
    \left\{(\frac{i}{10}, \frac{j}{10}, \frac{k}{10}) : i,j,k \in \{1,\dots,100\}\right\},
\end{equation*}
and each element represents
the triplet of the three constant functions $(f,\kappa,g)$.
Then, we randomly select 100,000 gird points out of 1M and solve the equation \eqref{eqn:Darcy} to obtain the corresponding solutions $p$ on $\my=541$ points.
The 100,000 data are split into 90,000 training data and 10,000 test data. 
Since the input functions are the triplet of constant functions, we use the corresponding grid as the input for DeepONets.
We employ the trunk and branch networks whose architectures are
$\vec{\bm{n}}_t = (2, 100, 100, 100, 200)$
and
$\vec{\bm{n}}_b = (3, 100, 100, 100, 201)$, respectively.
The hyperbolic tangent ($\tanh$) activation function is used for both
and the Xavier initialization scheme \cite{glorot2010understanding} and 
the \texttt{Adam} optimizer \cite{kingma2014adam} is utilized.

\begin{figure}[!h]
\centering
%\hfill
\begin{subfigure}[b]{0.49\textwidth}
         \centering
         \includegraphics[width=\textwidth]{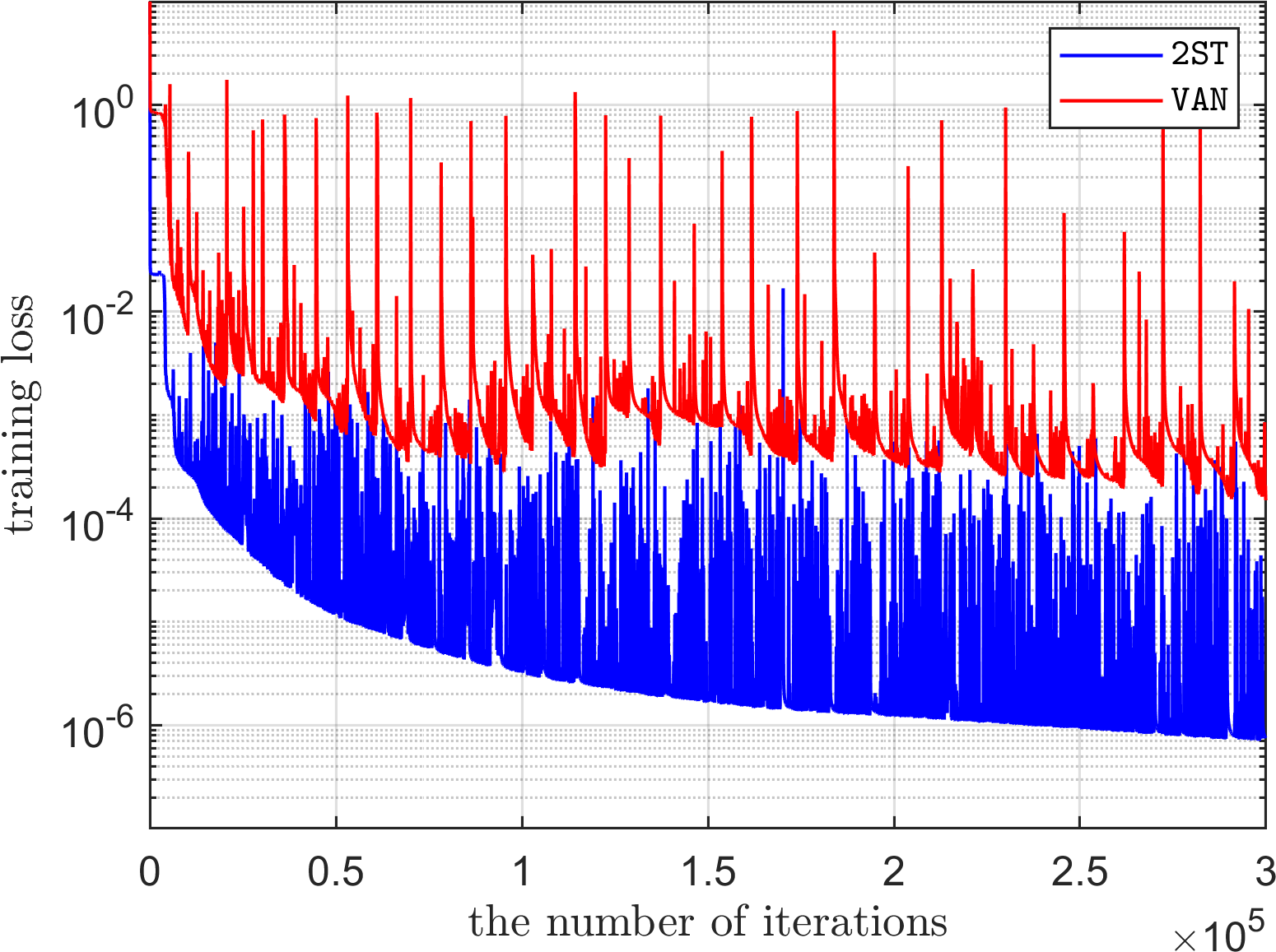}
         \caption{} \label{fig:ex3_loss_a}
\end{subfigure}
\begin{subfigure}[b]{0.49\textwidth}
         \centering
         \includegraphics[width=\textwidth]{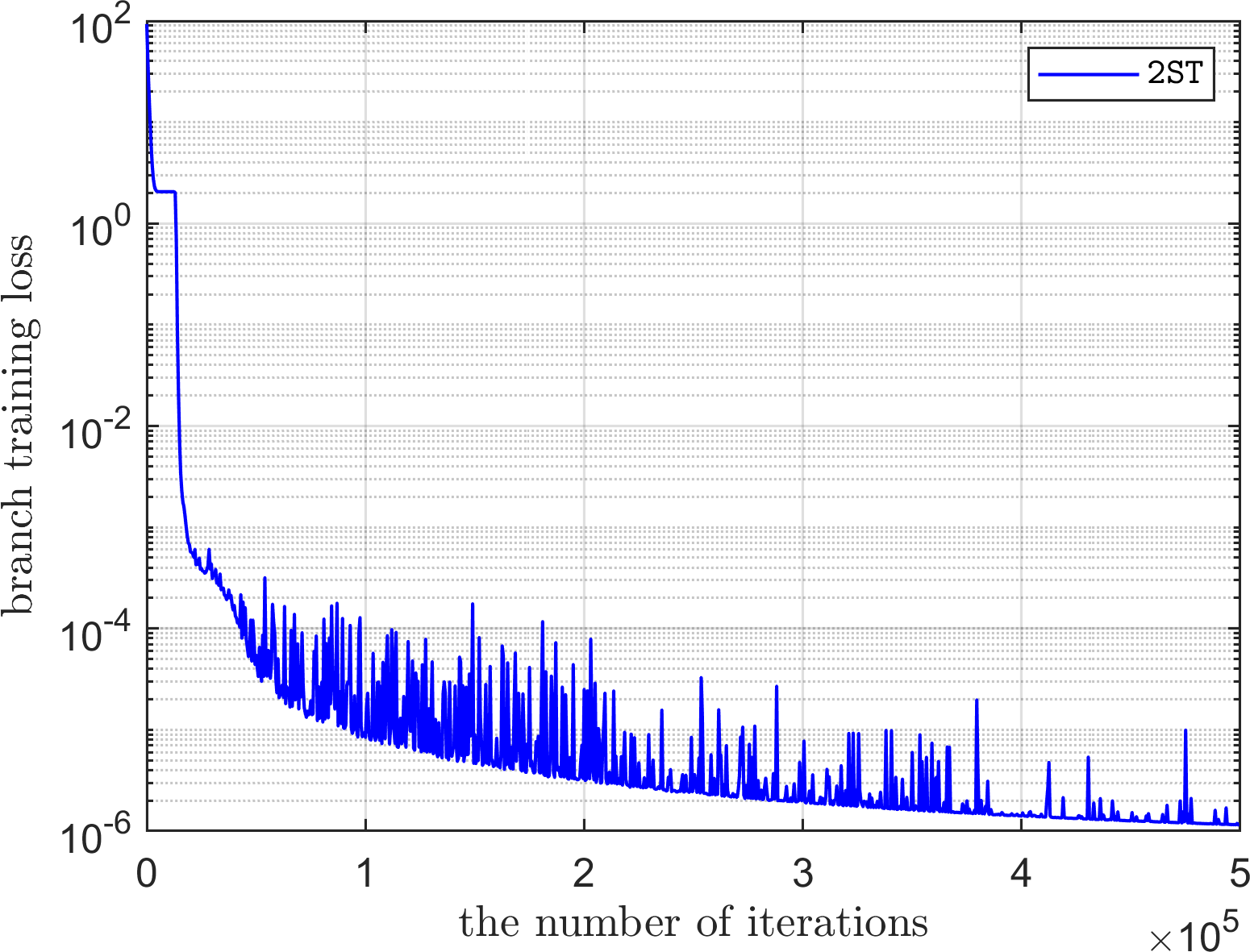}
         \caption{}%Approximated solution from ONet}
         \label{fig:ex3_loss_b}
\end{subfigure}
\caption{Example 5.3. (Left) The training loss versus the number of iterations by \texttt{VAN} and \texttt{2ST}. 
(Right) The branch training loss versus the number of iterations. This is the second step of the proposed two-step training method.}
\label{fig:ex3_loss}
\end{figure}

In Figure~\ref{fig:ex3_loss_a},
the training loss versus the number of iterations is reported.  Again, it is clearly seen that 
\texttt{2ST} can effectively minimize the trunk network loss reaching the level of $10^{-7}$ at the end of the training,
while \texttt{VAN} stagnates at the level of $10^{-4}$.
This again confirms that \texttt{2ST} effectively trains the trunk network when it is compared with \texttt{VAN}.
The remaining job for \texttt{2ST} is to train the branch network according to \eqref{train-branch}.

Figure~\ref{fig:ex3_loss_b} shows the branch loss in the second step of \texttt{2ST} with respect to the number of \texttt{Adam} iterations. 
Here, we use a learning rate scheduler that starts at $10^{-3}$ and reduces the learning rate by a factor of 2
for every 100K iteration.
It is observed that the branch loss is sufficiently minimized at the end of the training and reaches the level of $10^{-6}$.

\begin{figure}[htbp]
	\centering
\includegraphics[width=0.75\textwidth]{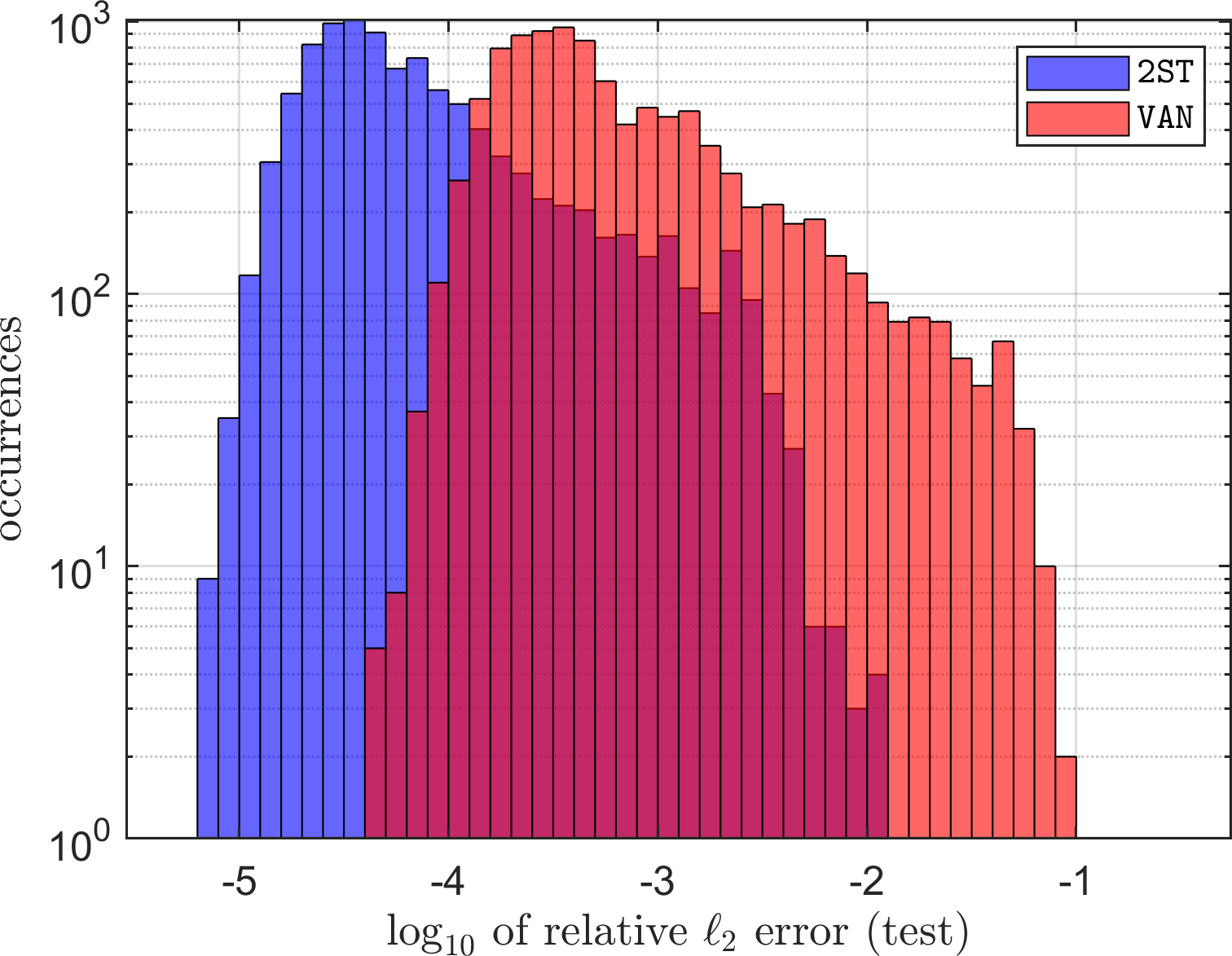}
	\caption{Example 5.3. 
        The histogram of $\log_{10}$ of the relative $\ell_2$ errors on the 10,000 test data obtained by 
        \texttt{2ST} (blue)
        and \texttt{VAN} (red).
	}
	\label{fig:ex3_test}
\end{figure}

Lastly, we report the histogram of $\log_{10}$ of the relative $\ell_2$ errors on the 10,000 test data
in Figure~\ref{fig:ex3_test}.
It is clearly observed that the DeepONet trained by \texttt{2ST} yields a much better generalization performance 
than the one by \texttt{VAN}.
The average relative test error by the proposed two-step method is 
$2.9\times 10^{-4}$,
while the one by the vanilla monolithic training is $2.5 \times 10^{-3}$.
This clearly indicates that how DeepONets are trained makes a significant impact on generalization performance. We emphasize that the only change we make is the training method,
while the network architecture, data, and initialization schemes were identical throughout.

\section{Conclusions}
In this study, we explored a novel training technique for DeepONets. The newly introduced sequential two-step training approach involves initial training of the trunk network that involves the Gram-Schmidt orthonormalization by means of QR-factorization, followed by the training of the branch network. The efficacy of the two-step training method was assessed through various numerical experiments, contrasting its performance against the conventional monolithic training approach, involving both forward and inverse Darcy problems within porous media contexts.
The efficacy and robustness of the proposed approach were clearly showcased in these representative examples, underscoring the significance of having robust training algorithms tailored to a specific neural network architecture.
Moreover, the significance of pre-training the trunk network was emphasized, as it provided valuable insights into the outcomes of the complete training process. This approach resulted in crucial improvements in accuracy while also simplifying the overall complexity of the training. 
Lastly, a generalization error estimate is established by leveraging the least-squares error analysis of \cite{cohen2013stability,cohen2019correction} in terms of the number of training data, the number of input and output sensor points, and the width of DeepONets.
% thereby expanding the range of applications for this approach.

\appendix

\section{Proof of Theorem~\ref{thm:trunk-interpolation}}
\label{app:them:trunk-interpolation}
\begin{proof}
    Let the rank of $\bm{U}$ be $r$
    and let $Z \Sigma_r V^\top$ be a SVD of $\bm{U}$
    where $Z \in \mathbb{R}^{\my \times r}$,
    $\Sigma_r \in \mathbb{R}^{r\times r}$
    and $V \in \mathbb{R}^{K \times r}$.
    Let
    \begin{align*}
        \tilde{A} = 
        \begin{cases}
            \begin{pmatrix}
            \Sigma_r V^\top \\
            \bm{0}_{(N-r)\times K}
            \end{pmatrix} & \text{if } N \ge r \\
            \Sigma_{1:N}V_{1:N}^\top  & \text{if } N < r 
        \end{cases},
    \end{align*}
    where
    $\Sigma_{1:s}$
    is a diagonal matrix of size $s\times s$ 
    obtained from $\Sigma_r$
    by collecting 
    the first $s$ rows and columns,
    and $V_{1:s}$
    is obtained by collecting the first $s$ columns of $V$.
    It then can be checked that 
    if the trunk network 
    satisfies
    \begin{equation} \label{eqn:thm:interpolation}
        \bm{\phi}_0^\top(y_i;\mu) = 
        \begin{cases}
            (Z^{(i)}, \bm{0}_{1\times (N-r)}) & \text{if } N \ge r \\
            Z^{(i)}_{1:N} & \text{if } N < r
        \end{cases},
    \end{equation}
    where $Z^{(i)}$ is the $i$-th row of $Z$ and $Z^{(i)}_{1:s}$ is the first $s$ entries of $Z^{(i)}$,
    the desired statement is obtained
    by letting $A^* = [\bm{0}, \tilde{A}]$
    as
    $\bm{\Phi}(\mu)A^* = \bm{\Phi}_0(\mu)\tilde{A} = Z_{1:\tilde{r}}\Sigma_{1:\tilde{r}}V_{1:\tilde{r}}^\top$
    where $\tilde{r} = \min\{N,r\}$
    and $\bm{\Phi}_0(\mu)$
    is the matrix whose $i$-th row is
    $\bm{\phi}_0^\top(y_i;\mu)$.
    
    For the rest of the proof, we explicitly construct a deep ReLU network satisfying \eqref{eqn:thm:interpolation}.
    We closely follow the construction that appeared in \cite{vardi2022on}.
    Note that for any distinct $y_1,\dots,y_{\my}$ in 
    $\mathbb{R}^{d_y}$,
    there exists a unit vector $v \in \mathbb{R}^{d_y}$ (see e.g. \cite{pmlr-v134-park21a})
    such that 
    \begin{align*}
        \sqrt{\frac{8}{\pi d}} \frac{1}{\my^2}
        \|y_i - y_j\|
        \le |v^\top (y_i - y_j)|
        \le \|y_i - y_j\|, \quad 
        \forall i \ne j.
    \end{align*}
    Let $\text{y}_i = \tilde{v}^\top y_i$
    where $\tilde{v} = 2\delta^{-1}\sqrt{\frac{\pi d}{8}} \my^2 v$,
    and 
    let $\text{y}_1 < \dots < \text{y}_{\my}$
    (after reordering if necessary).
    It then can be checked that $|\text{y}_i - \text{y}_j| \ge 2$ for all $i\ne j$.
    
    For $a < b$, let $N_{a,b}$ be
    a 3-layer ReLU network of width 2
    defined by      
    \begin{align*}
        N_{a,b}(\text{y}) = A^3 \sigma(A^2 \sigma(A^1 \text{y} + b^1) + b^2) + b^3,
    \end{align*}
    where $A^1 = \begin{pmatrix}
          -2 \\ 2
    \end{pmatrix}$,
    $b^1 = \begin{pmatrix}
        2a \\ -2b
    \end{pmatrix}$,
    $A^2 = -\begin{pmatrix}
        1 & 0 \\ 0 & 1
    \end{pmatrix}$,
    $b^2 = \begin{pmatrix}
        1 \\ 1 
    \end{pmatrix}$,
    $A^3 = \begin{pmatrix}
        1 & 1
    \end{pmatrix}$,
    $b^3 = -1$,
    which emulates the hat-like function, i.e.,
    \begin{align*}
        N_{a,b}(\text{y}) = \begin{cases}
            1 & \text{if } a \le \text{y} \le b \\
            0 & \text{if } \text{y} \le a - \frac{1}{2} \quad \text{or} \quad 
            b + \frac{1}{2} \le \text{y} \\
            2(\text{y} - (a-\frac{1}{2})) & \text{if } a - \frac{1}{2} < \text{y} < a \\
            -2(\text{y} - b) & \text{if } b < \text{y} < b + \frac{1}{2}
        \end{cases}.
    \end{align*}
    For $k=1,\dots,\tilde{r}$, let $N^{(k)}(\text{y}) = \sum_{j=1}^{\my} Z^{(j)}_k N_{\text{y}_j,\text{y}_{j+1}}(\text{y})$
    which satisfies 
    $N^{(k)}(\text{y}_j) = Z^{(j)}_k$
    for all $j \in [\my]$.
    The remaining task is to construct 
    a deep ReLU network $\bm{\phi}_0$ such that
    $\bm{\phi}_0(y;\mu) = (N^{(1)}(\tilde{v}^\top y),\dots, N^{(\tilde{r})}(\tilde{v}^\top y), \bm{0}_{N-\tilde{r}})^\top$.
    % \begin{align*}
    %     \bm{\phi}_0(y;\mu) = \begin{pmatrix}
    %         N^{(1)}(\tilde{v}^\top y) \\
    %         \vdots \\
    %         N^{(\tilde{r})}(\tilde{v}^\top y)
    %         \\
    %         \bm{0}_{N-\tilde{r}}
    %     \end{pmatrix}.
    % \end{align*}

    Let $\bm{z} = (z_1,\dots,z_{\tilde{r}})^\top$,
    $\bm{c} = (c_1,\dots,c_{\tilde{r}})^\top$,
    and consider 
    a 3-layer ReLU network $F_{a,b,\bm{c}}$ of width $\tilde{n}=2\tilde{r}+4$ defined by
    \begin{align*}
        F_{a,b,\bm{c}}(\begin{bmatrix}
            \text{y} \\ \bm{z}
        \end{bmatrix}) = \tilde{A}^3 \sigma(\tilde{A}^2 \sigma(\tilde{A}^1 \begin{bmatrix}
            \text{y} \\ \bm{z}
        \end{bmatrix} + \tilde{b}^1) + \tilde{b}^2) + \tilde{b}^3
        = \begin{pmatrix}
            \text{y} \\ \bm{z} + \bm{c}N_{a,b}(\text{y})
        \end{pmatrix},
    \end{align*}
    where $P = (1, -1)^\top$, 
    \begin{align*}
        \tilde{A}^1 = \begin{pmatrix}
        A^1 & \bm{0}_{2 \times \tilde{r}} \\
        P & \bm{0}_{2 \times \tilde{r}}  \\
        \bm{0}_{2\tilde{r} \times 1} & \textbf{P}_{\tilde{r}}
    \end{pmatrix} \in \mathbb{R}^{\tilde{n} \times ({\tilde{r}}+1)}
    \text{ with } 
    \textbf{P}_s = 
    \begin{pmatrix}
        P & \bm{0}_{2\times 1} & \cdots & \bm{0}_{2\times 1} \\
        \bm{0}_{2\times 1} & P & \cdots & \bm{0}_{2\times 1} \\
        \vdots & \vdots & \ddots & \vdots \\
        \bm{0}_{2\times 1} & \bm{0}_{2\times 1}
        & \cdots & P
    \end{pmatrix}
    \in \mathbb{R}^{2s \times s},
    \end{align*}
    $\tilde{b}^1 = \begin{pmatrix}
        b^1 \\ \bm{0}_{2(\tilde{r}+1)\times 1}
    \end{pmatrix}$,
    $\tilde{A}^2 = \begin{pmatrix}
        A^2 & \bm{0}_{2\times 2(\tilde{r}+1)}\\ 
        \bm{0}_{2(\tilde{r}+1) \times 2} & I_{2(\tilde{r}+1)}
    \end{pmatrix} \in \mathbb{R}^{\tilde{n} \times \tilde{n}}$,
    $\tilde{b}^2 = \begin{pmatrix}
        b^2 \\ \bm{0}_{2(\tilde{r}+1)\times 1}
    \end{pmatrix} \in \mathbb{R}^{\tilde{n}}$,
    $\tilde{A}^3 = \begin{pmatrix}
        \bm{0}_{1\times 2} & P^\top  & \bm{0}_{1\times 2\tilde{r}}
        \\ \bm{c}A^3 & \bm{0}_{\tilde{r} \times 2} & \textbf{P}_{\tilde{r}}^\top 
    \end{pmatrix} \in \mathbb{R}^{(\tilde{r}+1) \times \tilde{n}}$,
    $\tilde{b}^3 = \begin{pmatrix}
        0 \\ \bm{c}
    \end{pmatrix} \in \mathbb{R}^{\tilde{r}+1}$.

    Lastly, let us consider a 3-layer ReLU network
    $F^{(0)}_{a,b,\bm{c}}$ of width 4 defined by
    \begin{align*}
        F^{(0)}_{a,b,\bm{c}}(y) = \hat{A}^3 \sigma(\hat{A}^2 \sigma(\hat{A}^1 y + \hat{b}^1) + \hat{b}^2) + \hat{b}^3
        = \begin{pmatrix}
            x \\ \bm{c}N_{a,b}(\tilde{v}^\top y)
        \end{pmatrix},
    \end{align*}
    where
    $\hat{A}^1 = \begin{pmatrix}
        A^1 \\ P
    \end{pmatrix}\tilde{v}^\top \in \mathbb{R}^{4 \times d_y}$,
    $\hat{b}^1 = \begin{pmatrix}
        b^1 \\ \bm{0}_{2\times 1}
    \end{pmatrix}$,
    $\hat{A}^2=\begin{pmatrix}
        A^2 & \bm{0}_{2\times 2} \\
        \bm{0}_{2\times 2} & I_2
    \end{pmatrix}$,
    $\hat{b}^2 = \begin{pmatrix}
        b^2 \\ \bm{0}_{2\times 1}
    \end{pmatrix}$,
    $\hat{A}^3 = \begin{pmatrix}
        \bm{0}_{1\times 2} & P^\top \\
        \bm{c}A^3 & \bm{0}_{\tilde{r}\times 2}
    \end{pmatrix}$,
    $\hat{b}^3 = \begin{pmatrix}
        0 \\ \bm{c}
    \end{pmatrix}$.

    Let us consider 
    \begin{align*}
        \bm{\phi}_0(y;\mu^*):= 
        A
        F_{\text{y}_{\my-1},\text{y}_{\my}, Z_{1:\tilde{r}}^{(\my)}}
        \circ 
        \cdots 
        \circ 
        F_{\text{y}_2,\text{y}_{3}, Z_{1:\tilde{r}}^{(2)}}
        \circ
        F_{\text{y}_1,\text{y}_{2}, Z_{1:\tilde{r}}^{(1)}}^{(0)}(y),
    \end{align*}
    where 
    $A = \begin{bmatrix}
            \bm{0}_{\tilde{r}\times 1}
            & I_{\tilde{r}} \\
            \bm{0}_{(N-\tilde{r})\times 1}
            &  \bm{0}_{(N-\tilde{r})\times (N-\tilde{r})}
        \end{bmatrix}$.
    Then, it is a $(2\my+1)$-layer ReLU network 
    whose architecture is $\vec{\bm{n}}_t$
    as shown in the statement,
    and $\mu^*$ represents the corresponding 
    the network parameters.
    It then can be checked that 
    $\bm{\phi}_0(y;\mu^*) = (N^{(1)}(\tilde{v}^\top y),\dots, N^{(\tilde{r})}(\tilde{v}^\top y), \bm{0}_{N-\tilde{r}})^\top$,
    which completes the proof.
\end{proof}

\section{Proof of Theorem~\ref{thm:main}} \label{app:thm:main}

\begin{proof}
    Let $V_n(\mu) = \text{span}\{(\bm{\phi}(\cdot;\mu))_j : j=1,\dots,N+1\}$.
    % Suppose that $(\bm{\phi}(\cdot;\mu))_j$'s are linearly independent.
    The minimization problem of \eqref{train-trunk} is equivalent to  
    \begin{align*}
        \mu^* =\argmin_{\mu \in \text{M}} \| (\bm{\Phi}(\mu)\bm{\Phi}(\mu)^\dagger - I)\bm{U}\|_{2,2}^2,
    \end{align*}
    where $\text{M}$ is a feasible set for the trunk networks
    defined on Assumption~\ref{assmpt-trunk}.
    % It then follows from
    % \eqref{assumption-eqn-trunk-sampling}
    % that 
    % \begin{align*}
    %     \frac{1}{\my}\mathbb{E}\left[\|(\bm{\Phi}(\mu)\bm{\Phi}(\mu)^\dagger - I)\bm{u}\|^2\right]
    %     &= \mathbb{E}\left[ \bigg(u(Y) - \sum_{i=0}^N \hat{\phi}(Y;\mu;T_\mu) \hat{a}_j \bigg)^2\right]
    %     \\
    %     &= \|u - \Pi_{V_n(\mu)}u\|_{L^2_\omega}^2
    %     + \|
    % \end{align*}
    Since $\mu^* \in \text{M}$,
    there exists $T_{\mu^*}$
    such that 
    $\hat{\bm{\phi}}(\cdot;\mu^*,T_{\mu^*})$
    form
    orthonormal basis 
    in $L^2_\omega(\Omega_y)$,
    which we denote by 
    $\{\psi_{j}^*\}_{j=0}^{N}$.
    We then rewrite 
    the operator $\mathcal{G}$ of interest
    in terms of the new orthonormal basis, i.e.,
    \begin{align*}
        \mathcal{G}[f](y) = \sum_{j=0}^{\infty}
        \mathfrak{c}^*_j(f) \psi_j^*(y).
    \end{align*}
    Observe that
    \begin{align*}
        \|\mathcal{G}[f] - \tONet[f]\|_{L_\omega^2(\Omega_y)}^2
        &= 
        \|\mathcal{G}[f] - \mathcal{G}_N[f] + \mathcal{G}_N[f] - \tONet[f]\|_{L_\omega^2(\Omega_y)}^2 \\
        &= 
        \|\mathcal{G}[f] - \mathcal{G}_N[f]\|_{L_\omega^2(\Omega_y)}^2 + \|\mathcal{G}_N[f] - \tONet[f]\|_{L_\omega^2(\Omega_y)}^2 
        \\
        &=
        \|\mathcal{G}[f] - \mathcal{G}_N[f]\|_{L_\omega^2(\Omega_y)}^2
        + \sum_{j=0}^N (\mathfrak{c}^*_j(f) - \bm{c}_j(\bm{f};\theta^*))^2.
    \end{align*}
    The second term on the right-hand side of the above can be further bounded as follows:
    \begin{align*}
        &(\mathfrak{c}^*_j(f) - \bm{c}_j(\bm{f};\theta^*))^2
        \\
        &= (\mathfrak{c}^*_j(f) - \mathfrak{c}^*_j(f_k) + \mathfrak{c}^*_j(f_k) - \bm{c}_j(\bm{f};\theta^*) - \bm{c}_j(\bm{f}_k;\theta^*) + \bm{c}_j(\bm{f}_k;\theta^*))^2
        \\
        &\le 
        3\left\{(\mathfrak{c}^*_j(f) - \mathfrak{c}^*_j(f_k))^2 + (\mathfrak{c}^*_j(f_k) - \bm{c}_j(\bm{f}_k;\theta^*))^2 + (\bm{c}_j(\bm{f};\theta^*) - \bm{c}_j(\bm{f}_k;\theta^*))^2 \right\} 
        \\
        &\le 
        3\left\{L_j^2 d_{\mathcal{X}}^2(f,f_k) + (\mathfrak{c}^*_j(f_k) - \bm{c}_j(\bm{f}_k;\theta^*))^2 + L_{\sigma}^2L_{c}^2(K,N,\mx)\|\bm{f} -\bm{f}_k\|^2_{w,2} \right\}.
    \end{align*}
    Thus, we have
    \begin{align*}
        \|\mathcal{G}[f] - \tONet[f]\|_{L_\omega^2(\Omega_y)}^2
        &\le 
        \|\mathcal{G}[f] - \mathcal{G}_N[f]\|_{L_\omega^2(\Omega_y)}^2
        \\
        &\quad + 
        3\left\{L_\mathcal{G}^2 d_{\mathcal{X}}^2(f,f_k) + L_{\sigma}^2L_c^2 \|\bm{f} -\bm{f}_k\|^2_{w,2} \right\}
        \\
        &\quad +3\sum_{j=0}^N (\mathfrak{c}^*_j(f_k) - \bm{c}_j(\bm{f}_k;\theta^*))^2.
    \end{align*}
    Recall that since the trunk networks are an orthogonal basis,
    the optimal solution $A^*$ of \eqref{train-trunk} is the least squares solution as $p=2$.
    Under Assumptions~\ref{assmpt-operators} and~\ref{assmpt-trunk},
    it follows from Theorem 2 of \cite{cohen2013stability,cohen2019correction}
    that we have
    \begin{align*}
        \mathbb{E}\left[\sum_{j=0}^N (\mathfrak{c}^*_j(f_k) - \bm{c}_j(\bm{f}_k;\theta^*))^2\right]
        \lesssim C'(\my,r_t)\|\mathcal{G}[f_k] - \mathcal{G}_N[f_k]\|_{L_\omega^2(\Omega_y)}^2
        + \my^{-r_t},
    \end{align*}
    where $C'(m_y,r_t) = \frac{6\log (3/2) - 2}{(1+r_{t})\log \my}$. 
    % \shlee{What is C?}
    By combining the above with 
    Assumptions~\ref{assmpt-operators} and~\ref{assmpt-inputspace}, we have
    \begin{align*}
        \mathbb{E}\left[\|\mathcal{G}[f] - \tONet[f]\|_{L_\omega^2}^2\right]
        \lesssim
        C(\my,r_t)N^{-r_{\mathcal{G,X},\mu^*}}
        + \my^{-r_t} +
        K^{-\alpha} + \mx^{-s},
    \end{align*}
    which completes the proof.
\end{proof}

% Since $\{y_j\}$'s are i.i.d. samples from $\omega$,
%     we have 
%     \begin{align*}
%         \sum_{k=1}^K \inf_{g \in V_n(\mu)} \|g - u_k\|_{L^2_\omega}^2
%         =
%         \min_{A}
%         \frac{1}{K\my}\mathbb{E}\left[\|\bm{\Phi}(\mu)A - \bm{U}\|_{2,2}^2\right]
%         \le \frac{1}{K\my}\mathbb{E}\left[\|\bm{\Phi}(\mu)(\bm{\Phi}(\mu)^\dagger \bm{U}) - \bm{U}\|_{2,2}^2\right]
%     \end{align*}
%     where the expectation is taken over the i.i.d. samples.
%     It then follows from \cite{cohen2013stability,cohen2019correction} that for each $\mu \in \text{M}$,
%     \begin{align*}
%         \frac{1}{K}\sum_{k=1}^K \mathbb{E}\bigg[ \|\tONet[\bm{f}_k;\mu] - u_k \|_{L^2_\omega}^2\bigg]
%         &\lesssim
%         (1 + e(\my)) \frac{1}{K}\sum_{k=1}^K \inf_{g \in V_n(\mu)} \| g - u_k\|_{L^2_\omega}^2 + \my^{-r}.
%         % \\
%         % &\lesssim
%         % \frac{1}{K\my}\|(\bm{\Phi}(\mu)\bm{\Phi}(\mu)^\dagger - I)\bm{U}\|_{2,2}^2
%         % + \my^{-r}.
%     \end{align*}
% \input appendix

\bibliographystyle{siamplain}
\bibliography{refs}

\end{document}